\documentclass[12pt]{amsart}

\usepackage[scale=.86]{geometry}
\usepackage{geometry} 
\usepackage{graphicx}
\usepackage{bbm}
\usepackage{mathrsfs}
\theoremstyle{plain}

\usepackage{CJK}
\usepackage{amsfonts,amssymb,comment,amsmath,mathrsfs,multirow,booktabs}

\usepackage{hyperref}  

\hypersetup{
	colorlinks = true,  
	linkcolor = blue,  
	citecolor = blue,   
	urlcolor  = red     
}

\usepackage{color,latexsym,amsfonts}
 \numberwithin{equation}{section}

 \newtheorem{thm}{Theorem}[section]
 \newtheorem{cor}[thm]{Corollary}
 \newtheorem{lem}[thm]{Lemma}
 \newtheorem{prop}[thm]{Proposition}
 \theoremstyle{remark}
 \newtheorem{rem}[thm]{Remark}
 \theoremstyle{example}

 \def\ar{\!\!\!&}

 \def\mcr{\mathscr}\def\mbb{\mathbb}\def\mbf{\mathbf}

 \def\beqlb{\begin{eqnarray}}\def\eeqlb{\end{eqnarray}}
 \def\beqnn{\begin{eqnarray*}}\def\eeqnn{\end{eqnarray*}}
 \def\d{{\mbox{\rm d}}}

 \def\eqref#1{{\rm(\ref{#1})}}

\begin{document}



\centerline{\Large\textbf{ Threshold Diffusions}}

\bigskip

\centerline{Lina Ji\,\footnote{Supported by NSFC grant (No. 12301167), Guangdong Basic and Applied Basic Research Foundation (No. 2022A1515110986) and Shenzhen National Science Foundation (No. 20231128093607001).}, Chuyang Li, Xiaowen Zhou\,\footnote{Supported by Natural Sciences and Engineering Research Council of Canada (RGPIN-2021-04100).}}

\medskip

\centerline{\it MSU-BIT-SMBU Joint Research Center of Applied Mathematics, }

\centerline{\it Shenzhen MSU-BIT University, Shenzhen 518172, P.R. China.}
\centerline{\it  School of Mathematics and Systems Science, }
\centerline{\it Guangdong Polytechnic Normal University, Guangzhou 510665, P.R. China.}
 
\centerline{\it Department of Mathematics and Statistics, Concordia University}

\centerline{\it 1455 De Maisonneuve Blvd. W., Montreal, Canada.}

\centerline{E-mails: \tt jiln@smbu.edu.cn, lichuyang@gpnu.edu.cn,}
\centerline{\tt xiaowen.zhou@concordia.ca}

\bigskip

{\narrower{\narrower
		
		\centerline{\textbf{Abstract}}
		
	{ 	We propose threshold diffusion processes as unique solutions to stochastic differential equations with step-function coefficients, and obtain explicit expressions for the conditional Laplace transform of the hitting times and the potential measures. Applying these results, we further discuss their asymptotic behaviours such as the stationary distributions and the  escape probabilities. }

		\medskip
		
		
		\par}\par}

\bigskip

\section{Introduction and Main Results}

\subsection{Introduction}
Let $(\Omega, \mcr{F}, (\mcr{F}_t)_{t \ge 0}, \mbf{P})$ be a filtered probability space satisfying the usual conditions, and $B := (B_t)_{t \ge 0}$ be a one-dimensional Brownian motion adapted to $(\mcr{F}_t)_{t \ge 0}$. A {\it multi-regime threshold diffusion process} $X := (X_t)_{t \ge 0}$ is defined as the unique strong solution to the following threshold stochastic differential equation:
\beqlb\label{Xn}
\d X_t \ar=\ar b(X_t)\d t + \sigma(X_t)\d B_t 
\eeqlb
with step-function coefficients 
\beqlb\label{b}
\begin{cases}
b(x) = \mu_0{\bf 1}_{(-\infty, a_1]}(x) + \sum_{i = 1}^{n-1}\mu_i {\bf 1}_{(a_i, a_{i+1}]}(x) + \mu_{n}{\bf 1}_{(a_n, \infty)}(x), \cr
\sigma(x) = \sigma_0{\bf 1}_{(-\infty, a_1]}(x) + \sum_{i = 1}^{n-1}\sigma_i {\bf 1}_{(a_i, a_{i+1}]}(x) + \sigma_{n}{\bf 1}_{(a_n, \infty)}(x),
\end{cases}
\eeqlb
where $\mu_i \in \mbb{R},\ i = 0, 1, \cdots, n$ are the drift parameters, $\sigma_i \ge 0,\ i = 0, 1, \cdots, n$ are the diffusion parameters, the  threshold vector $(a_1, \cdots, a_n)$ satisfying  $-\infty< a_1  <\cdots< a_n < \infty$ is a partition of $\mbb{R}$ and ${\bf 1}_{A}(x) = {\bf 1}(x \in A)$   denotes the indicator function. This is a typical diffusion process with discontinuous coefficients. Without loss of generality, we assume that $\sigma_i >0$ for any $i = 0, 1, \cdots, n$. The existence and strong uniqueness of the solution $X$ to \eqref{Xn}  can be established using the splicing method; also see \'{E}tor\'{e} and Martinez \cite{EM18} and Le Gall \cite{LG84}, which establish that \eqref{Xn} admits a unique strong solution. It is well known that continuous functions can be approximated by step functions. Thus, we consider diffusion processes with step-function coefficients in this paper, which can approximate diffusion processes with continuous coefficients, see Lejay and Martinez \cite{LM06}.

Many applications can be found for such diffusion processes with piecewise constant coefficients. It can be seen as an alternative to the model studied in Mota and Esqu\'{i}vel \cite{ME14},
which is a continuous-time version of the self-exciting threshold autoregressive models, a subclass of the discrete-time threshold autoregressive models presented in Tong \cite{T83,T11}. For the  single-threshold  ($n = 1$) diffusion, Lejay and Pigato \cite{LP20} established  the maximum likelihood estimate for its drift, and Chen et al. \cite{CWZ24} derived an explicit expression its transition density. This enables direct computation of value functions for control problems akin to those in \cite{M84}, with recent extensions rigorously explored in Wang et al. \cite{WXYYZ25}. See also Zhao and Xi \cite{ZX21} and references therein for related work.


Parallel developments have emerged in related work, notably for  threshold Ornstein-Uhlenbeck (TOU) processes. 
Hu and Xi \cite{HX22} studied the parameter estimation of TOU processes, with further extensions by Han and Zhang \cite{HZ23a, HZ23b} and Zhang \cite{Z24}. Further applications of threshold processes in financial contexts can be found in the works of Su and Chan \cite{SC15}, Ferrari and Vargiolu \cite{FV20}, Li et al. \cite{LMH23}, and Bernis et al. \cite{BBSS21}, and references therein.

The classical theory of one-dimensional diffusions with continuous coefficients is well-established; see, e.g., Itô and McKean \cite{IM96}, Karatzas and Shreve \cite{KS91}, and Rogers and Williams \cite{RW00}. These results include closed-form expressions for the speed measure and scale function, appearing in works such as Borodin and Salminen \cite[page 17]{BS02}, Rogers and Williams \cite[page 270]{RW00}, Karlin and Howard \cite[pages 194-195]{KH81}. It is worth mentioning that the stationary measures of classical diffusion processes can be expressed in terms of the speed measure, which it related to the scale function. This formulation typically requires continuous coefficients. In this paper, we prove that this result also holds for the discontinuous-coefficient model \eqref{Xn} defined by \eqref{b}.

This paper concerns the properties 
such as the conditional Laplace transform of hitting times and the potential measure for the above-mentioned process $X$:
\beqlb\label{e_q}
\mbf{P}_x\{X_{e_q} \in \d z\} = \int_0^\infty q e^{-qt}\mbf{P}_x\{X_t \in \d z\}\d t,
\eeqlb
where $e_q$ is an exponential random variable with rate $q > 0$, independent of $X$. Then the stationary distribution (when $\mu_0 > 0$ and $\mu_n < 0$) is  obtained by taking a limit in \eqref{e_q} as $q \rightarrow 0.$  Moreover, the escape probabilities for $X$ to $-\infty$ and to $+\infty$   (when $\mu_0 < 0$ and $\mu_n > 0$) are also obtained.

\subsection{Main Results} 
Prior to presenting the main results, we introduce the following notation. For $q\geq 0$ and $i = 0, 1, \cdots, n$, define 
\beqlb\label{delta}
 l_i := \frac{\sqrt{2q\sigma_i^2+ \mu_i^2} }{\sigma_i^2}, \qquad \delta_i^- := l_i - \frac{\mu_i}{\sigma_i^2}, \qquad \delta_i^+ := l_i +  \frac{\mu_i}{\sigma_i^2}.
\eeqlb
For constants $b_i^+, c^+_{i}, b_i^-$ and $c^-_{i}$ to be determined in the next section, let
\beqlb\label{g_q^+}
g_q^+(x) \ar=\ar e^{\delta_0^-(x-a_1)}\mbf{1}_{(-\infty, a_1]}(x)  + \sum_{i=1}^{n-1} b_i^+ \bigg( (1-c^+_{i})e^{\delta_{i}^-(x-a_i)}+c^+_{i}e^{-\delta_{i}^+(x-a_i)}\bigg)\mbf{1}_{(a_{i}, a_{i+1}]}(x)\cr
\ar\ar + b_n^+\bigg( (1-c^+_{n})e^{\delta_{n}^-(x-a_n)}+c^+_{n}e^{-\delta_{n}^+(x-a_n)}\bigg)\mbf{1}_{(a_n, \infty)}(x) 
\eeqlb
and
\beqlb\label{g_q^-}
g_q^-(x) \ar=\ar e^{-\delta_{n}^+(x-a_n)}\mbf{1}_{(a_n, \infty)}(x) + \sum_{i=2}^{n}b_i^-\bigg( c^-_{i}e^{\delta_{i-1}^-(x-a_i) }+(1-c^-_{i})e^{-\delta_{i-1}^+(x-a_i)}\bigg)\mbf{1}_{(a_{i-1}, a_{i}]}(x)\cr
\ar\ar + b_1^-\bigg( c^-_{1}e^{\delta_0^-(x-a_1) }+(1-c^-_{1})e^{-\delta_0^+(x-a_1)}\bigg)\mbf{1}_{(-\infty, a_1]}(x)
\eeqlb
 satisfying $g_q^{\pm} \in C^1(\mbb{R})\cap C^2(\mbb{R}\setminus \cup_{i = 1}^n\{a_i\})$. 

\begin{lem}\label{l0525}
	The processes $\left(e^{-qt}g_q^\pm(X_t) \right)_{t \ge 0}$ are $(\mcr{F}_t)_{t \ge 0}$ local martingales.
\end{lem}
\proof
It is easy to check that $g_q^{\pm}$ are the fundamental solutions to
$\mathcal{L}g(x) = qg(x)$, (see \cite[page 128]{IM96}), where
\beqlb\label{L}
\mathcal{L}g(x) \ar=\ar b(x)g'(x) + \frac{1}{2} \sigma(x)^2g''(x).
\eeqlb
Here $b(x)$ and $\sigma(x)$ are given by \eqref{b}. By Meyer-It\^{o}'s formula (see, e.g., Protter \cite[Theorem 70 of Chapter IV]{P04}) we have
\beqnn
e^{-qt}g(X_t) \ar=\ar g(X_0) + \int_0^t e^{-qs}\left[\mathcal{L}g(X_s) - qg(X_s)\right]\d s + local\ mart. 
\eeqnn 
Then the result follows.
\qed 
 
We now derive the conditional Laplace transform of the first passage times for $X$. Letting $\tau_y := \inf\{t \ge 0: X_t = y\}$ for $y\in\mathbb{R}$, we have the following result.

\begin{thm}
	For any $y \le x \le z,$ we have
	\beqlb\label{eq0305}
	\mbf{E}_x\left[e^{-q\tau_y}{\bf 1}_{\{\tau_y < \tau_z\}}\right] \ar=\ar \frac{g_q^-(x)g_q^+(z) - g_q^+(x)g_q^-(z)}{g_q^-(y)g_q^+(z) - g_q^+(y)g_q^-(z)}
	\eeqlb
	and
	\beqlb\label{eq0305a}
	\mbf{E}_x\left[e^{-q\tau_z}{\bf 1}_{\{\tau_z < \tau_y\}}\right] = \frac{g_q^-(x)g_q^+(y) - g_q^+(x)q_q^-(y)}{g_q^-(z)g_q^+(y) - g_q^+(z)g_q^-(y)}.
	\eeqlb
\end{thm}
\proof 
For $y \le x \le z,$ by Lemma \ref{l0525} it holds that
\beqnn
g_q^-(x) \ar=\ar \mbf{E}_x[e^{-q(\tau_y\wedge\tau_z)}g_q^-(X_{\tau_y\wedge\tau_z})]\cr
\ar=\ar  g_q^-(y)\mbf{E}_x[e^{-q\tau_y}\textbf{1}_{\{\tau_y < \tau_z\}}] + g_q^-(z)\mbf{E}_x[e^{-q\tau_z}\textbf{1}_{\{\tau_z < \tau_y\}}]
\eeqnn
and
\beqnn
g_q^+(x) \ar=\ar \mbf{E}_x[e^{-q(\tau_y\wedge\tau_z)}g_q^+(X_{\tau_y\wedge\tau_z})]\cr
\ar=\ar  g_q^+(y)\mbf{E}_x[e^{-q\tau_y}\textbf{1}_{\{\tau_y < \tau_z\}}] + g_q^+(z)\mbf{E}_x[e^{-q\tau_z}\textbf{1}_{\{\tau_z < \tau_y\}}].
\eeqnn
The result follows by solving the above two equations.
\qed

Notice that $\lim_{x \rightarrow \infty}g_q^-(x) = 0$ and $\lim_{x \rightarrow-\infty} g_q^+(x) = 0.$ By letting $z \rightarrow \infty$ in \eqref{eq0305} and letting $y \rightarrow -\infty$ in \eqref{eq0305a}, we get the following result.

\begin{cor}\label{exptau}
	For any $y \le x \le z,$ we have
	\beqnn
	\mbf{E}_x\left[e^{-q\tau_y}\right] = \frac{g_q^-(x)}{g_q^-(y)}
	\qquad \text{and} \qquad
	\mbf{E}_x \left[e^{-q\tau_z}\right] = \frac{g_q^+(x)}{g_q^+(z)}.
	\eeqnn
\end{cor}

In what follows, we study the scale function for \eqref{Xn}. Letting $q \rightarrow 0$ in \eqref{eq0305}, we  obtain the following result. However, since such an approach is computationally intensive, we will instead provide a proof using probabilistic methods later.

 The following result is about the potential measure of $X$.
\begin{thm}\label{tpm}

	\begin{itemize}
		\item[(i)] For any $z \le a_1$, we have 
	\beqnn
	\mbf{P}_{x}\{X_{e_q}\in \d z\}= \frac{q e^{-\delta_0^+(a_1 - z)}}{l_0\sigma_0^2 b_1^- (1-c^-_1)} \begin{cases}
		g_q^-(x) \d z, &  x\ge z,\cr
		g_q^-(z) e^{-\delta_0^-(z - x)} \d z, &    z\ge x;
	\end{cases}
	\eeqnn
\item [(ii)]
	For any $z \in (a_i, a_{i+1}]$ with $i = 1, \cdots, n-1$, we have
	\beqnn
	\mbf{P}_x\{X_{e_q}\in \d z\} = \frac{qe^{-2\mu_i(a_{i+1} - z)/\sigma_i^2}}{C_i l_i\sigma_i^2g_q^+(a_i)g_q^-(a_{i+1})}\begin{cases}
		g_q^+(x)g_q^-(z) \d z, & x \le z,\cr
		g_q^-(x)g_q^+(z)\d z,   & x > z,
	\end{cases} 
	\eeqnn
	where
	\beqlb\label{Ci}
	C_i := (1-c^+_{i})(1- c^-_{i+1}) e^{\delta^-_i(a_{i + 1} - a_i)} - c^+_{i}c^-_{i+1} e^{- \delta^+_{i}(a_{i+1} - a_{i})};
	\eeqlb 
	\item[(iii)] For any $z > a_n$, we have
	\beqnn
	\mbf{P}_{x}\{X_{e_q}\in \d z\}
	\ar=\ar 	\frac{qe^{-\delta_n^-(z - a_n)}}{l_n\sigma_n^2 b_n^+(1 - c_n^+)} \begin{cases}
		g_q^+(x) \d z, & z \ge x,\\
		g_q^+(z) e^{-\delta^+_n(x - z)} \d z, & x \ge z.
	\end{cases}
	\eeqnn	
	\end{itemize}
\end{thm}


 For the case of $\mu_0 > 0$ and $\mu_n < 0$, the process $X$ admits a stationary distribution, which can be obtained by letting $q \rightarrow 0$ in Theorem \ref{tpm}.

The following result presents the closed-form expression for the stationary distribution of process $X$.
\begin{cor}\label{main1}
	Assuming that $\mu_0 > 0$ and $\mu_n < 0$, for any $x \in \mathbb{R}$ we have
	\beqnn
	\ar\ar \lim_{t \rightarrow \infty}\mbf{P}_x\{X_t \in \d z\}\cr
	\ar\ar\quad = \begin{cases}
		\frac{2}{\sigma_0^2\overline{F}_{1}}e^{-2\mu_0(a_1 - z)/ \sigma_0^2}\d z, & z \le a_1,\\
		\frac{2}{\sigma_i^2\overline{F}_{1}}\exp\left\{\sum_{\ell = 1}^{i-1}\frac{2\mu_{\ell}(a_{\ell + 1} - a_{\ell})}{\sigma_{\ell}^2} + \frac{2\mu_i(z - a_i)}{\sigma_i^2}\right\}\d z, & z \in (a_i, a_{i+1}], i = 1, \cdots, n,
	\end{cases}
	\eeqnn
	where $\sum_{\ell = 1}^0 = 0$, $a_{n + 1} :=+\infty$ and
 \beqlb\label{overlineF1}
	 \overline{F}_1 \ar=\ar \sum_{k = 1}^n \exp\left\{\sum_{j = 1}^{k-1}\frac{2\mu_j(a_{j+1} - a_j)}{\sigma_j^2}\right\}\cr
	\ar\ar \times \left(\frac{1}{\mu_{k-1}}{\bf 1}_{\{\mu_{k-1} \neq 0\}} - \frac{1}{\mu_k}{\bf 1}_{\{\mu_k \neq 0\}} + \frac{2(a_{k+1} - a_k)}{\sigma_k^2}{\bf 1}_{\{\mu_k = 0\}}\right).
	\eeqlb 
\end{cor}

\begin{rem}
	(i) For the case of $n=1$ in \eqref{Xn} with $\mu_0 > 0$ and $\mu_1 < 0$, the stationary distribution of $X$ is given as follows: for any $x \in \mbb{R}$,
	\beqnn
	\lim_{t \rightarrow \infty}\mbf{P}_x\{X_t \in \d z\} = \begin{cases}
		\frac{2e^{2\mu_0(z-a_1)/ \sigma_0^2}}{\sigma_0^2 \left(\frac{1}{\mu_0} - \frac{1}{\mu_1}\right)}\d z, & z \le a_1,\\	
		\frac{2e^{2\mu_1(z-a_1)/ {\sigma_1}^2}}{\sigma_1^2\left(\frac{1}{\mu_0} - \frac{1}{\mu_1}\right)}\d z, & z > a_1.
	\end{cases} 
	\eeqnn 
	The result recovers \cite[Proposition 4.7]{CWZ24};

	(ii) The scale function $\phi(\cdot)$ for $X$ is defined as
		\beqlb\label{psi}
	\ar\ar \phi(x) :=\int_{a_1}^x \exp\left\{-\int_{a_1}^y\frac{2b(z)}{\sigma(z)^2}\d z\right\}\d y
	\eeqlb 
	with the convention $\int_a^b = -\int_b^a$ if $a > b,$ which is a solution to $\mathcal{L}\phi(x) = 0$, see, e.g., \cite[Definition 46.10]{RW00}.  Then 
	the associated speed measure $m(x)\d x$ can be defined as 
	\beqnn
	m(x)\d x \ar:=\ar \frac{\d x}{\sigma(x)^2\phi'(x)}\cr
	 \ar=\ar \begin{cases}
		\frac{1}{\sigma_0^2 } e^{2\mu_0(x-a_1)/ \sigma_0^2}\d x, & x \le a_1,\\
		\frac{1}{\sigma_i^2}\exp\left\{\int_{a_1}^x 2b(z)/\sigma(z)^2 d z\right\}\d x,  & x \in (a_i, a_{i+1}], i = 1, \cdots, n, 
	\end{cases}\cr
\ar=\ar \begin{cases}
	\frac{1}{\sigma_0^2 } e^{2\mu_0(x-a_1)/ \sigma_0^2}\d x, & x \le a_1,\\
	\frac{1}{\sigma_i^2}\exp\left\{\sum_{\ell = 1}^{i-1}\frac{2\mu_\ell(a_{\ell + 1} - a_{\ell})}{\sigma_\ell^2} + \frac{2\mu_i(x - a_i)}{\sigma_i^2}\right\}\d x,  & x \in (a_i, a_{i+1}], i = 1, \cdots, n 
\end{cases}
	\eeqnn
	with $a_{n + 1} = +\infty.$ 
Note that $[\int_{\mbb{R}}m(x)\d x]^{-1} m(x)\d x$ coincides with the stationary distribution in Corollary \ref{main1}. In fact, the stationary distribution for a diffusion process with continuous coefficients is given by $m(x)\d x$, defined as above, 
see, e.g., \cite[Theorems 52.1 and 53.1]{RW00} and \cite[pages 194-195]{KH81}. This result also holds for our model with discontinuous coefficients.
\end{rem}

Under the condition of $\mu_0 < 0$ and $\mu_n > 0$, the process $X$ does not admit a stationary distribution. Since $t \rightarrow \infty$, $X$ almost surely converges to $+ \infty$ or $- \infty$. In the following we derive the probability of convergence to $-\infty$, which does not depend on the initial value. We first consider a sequence $\{(A_i, B_i): i = 1, \cdots, n\}$ satisfying the following conditions:
\beqlb\label{ABi}
\begin{cases}
	A_i =  \frac{\mu_n/\sigma_n^2}{\mu_i/\sigma_i^2{\bf 1}_{\{\mu_i \neq 0\}} + {\bf 1}_{\{\mu_i = 0\}}} e^{\sum_{\ell = i}^{n-1}\mu_\ell (a_{\ell + 1} - a_\ell)/\sigma_{\ell}^2}, & i = 1, \cdots, n-1,\\ 
	B_i = \Big(A_{i+1} + B_{i + 1} - \frac{\mu_{i+1}/\sigma_{i+1}^2}{\mu_i/\sigma_i^2{\bf 1}_{\{\mu_i \neq 0\}} + {\bf 1}_{\{\mu_i = 0\}}} A_{i+1}& \\
	\qquad\quad + 2\frac{\mu_{i+1}}{\sigma_{i+1}^2}(a_{i+1} - a_i)A_{i+1}{\bf 1}_{\{\mu_i = 0\}} \Big)e^{-\mu_i(a_{i+1} - a_i)/\sigma_i^2}, & i = 1, \cdots, n-1,\\ 
	A_n = 1, \quad B_n = 0.
\end{cases}
\eeqlb

The following result establishes the probability of escape to $-\infty.$
\begin{cor}\label{t0227}
	If  $\mu_0 < 0$ and $\mu_n > 0$, then for any initial value $y \in \mbb{R}$, we have
	\beqnn
	\ar\ar\mbf{P}_y\left\{\lim_{t \rightarrow \infty} X_t = -\infty\right\}\cr
	\ar\ar\quad = \begin{cases}
		\frac{ A_{i}(y) + B_{i}(y)}{\left(1 - \frac{\mu_1/\sigma_1^2{\bf 1}_{\{\mu_1 \neq 0\}} + {\bf 1}_{\{\mu_1 = 0\}}}{\mu_0/\sigma_0^2}\right)A_1 + B_1}e^{\frac{\mu_{i-1}(a_i - y)}{\sigma_{i-1}^2} -\sum_{\ell = 1}^{i-1}\frac{\mu_\ell(a_{\ell + 1} - a_\ell)}{\sigma_\ell^2}}, & y \in (a_{i-1}, a_i],\\
		& i = 1, \cdots, n,\\
		\frac{1}{\left(1 - \frac{\mu_1/\sigma_1^2{\bf 1}_{\{\mu_1 \neq 0\}} + {\bf 1}_{\{\mu_1 = 0\}}}{\mu_0/\sigma_0^2}\right)A_1 + B_1}e^{-\sum_{\ell = 1}^{n-1}\frac{\mu_\ell(a_{\ell + 1} - a_\ell)}{\sigma_\ell^2}-\frac{2\mu_n(y - a_n)}{\sigma_n^2}}, & y > a_n 
	\end{cases}
	\eeqnn 
	with $a_0 := -\infty$ and $\sum_{\ell = 1}^0 := 0$, where $(A_1, B_1)$ are given by \eqref{ABi} and for $y \in (a_{i-1}, a_i]$ with $i = 1, \cdots, n$,
 \beqnn
	\begin{cases}
		A_{i}(y) = \frac{\mu_i/\sigma_i^2 {\bf 1}_{\{\mu_i \neq 0\}} + {\bf 1}_{\{\mu_i = 0\}}}{\mu_{i-1}/\sigma_{i-1}^2} A_i e^{\frac{\mu_{i-1}(a_i - y)}{\sigma_{i-1}^2}}{\bf 1}_{\{\mu_{i-1} \neq 0\}}, \\ 
		B_{i}(y) = \left(A_i + B_i -\frac{\mu_i/\sigma_i^2 {\bf 1}_{\{\mu_i \neq 0\}} + {\bf 1}_{\{\mu_i = 0\}}}{\mu_{i-1}/\sigma_{i-1}^2} A_i{\bf 1}_{\{\mu_{i-1} \neq 0\}} + \frac{2\mu_i(a_i - y)}{\sigma_i^2}A_i {\bf 1}_{\{\mu_{i-1} = 0\}}\right)e^{-\frac{\mu_{i-1}(a_i - y)}{\sigma_{i-1}^2}}.
	\end{cases}
	\eeqnn 
\end{cor}

\begin{rem}
	(i) The probability $\mbf{P}_y\{\lim_{t \rightarrow \infty} X_t = + \infty\}$ can be obtained directly by Corollary \ref{t0227}, since 
	\beqnn
	\mbf{P}_y\left\{\lim_{t \rightarrow \infty} X_t = + \infty\right\} + \mbf{P}_y\left\{\lim_{t \rightarrow \infty} X_t = - \infty\right\} = 1;
	\eeqnn
	(ii) For the case of $n = 1$ in \eqref{Xn} with $\mu_0 < 0$ and $\mu_1 > 0$, for any $y \in \mbb{R}$ we have 
	\beqnn
	\mbf{P}_y\left\{\lim_{t \rightarrow \infty}X_t = -\infty\right\} = \begin{cases}
		1 + \frac{ \mu_1/\sigma_1^2 e^{2\mu_{0}(a_1 - y)/\sigma_{0}^2}}{ \mu_0/\sigma_0^2 - \mu_1/\sigma_1^2}, & y \le a_1,\\
		\frac{\mu_0/\sigma_0^2}{\mu_0/\sigma_0^2 - \mu_1/\sigma_1^2 }e^{ -2\mu_1(y - a_1)/\sigma_1^2}, & y > a_1.
	\end{cases}
	\eeqnn 
	The result is identical to \cite[Subsection 6.5]{LP20} when $y = a_1 = 0$.
\end{rem}

The remainder of this paper is organized as follows. Some basic results are presented in Section \ref{s}. The proof of Theorem \ref{tpm} is given in Section \ref{S3}, which characterizes the distribution of $X_{e_q}$. The proof of Corollary \ref{main1} is provided in Section \ref{SI}, while Section \ref{SP} contains the proof of Corollary \ref{t0227}.

\section{Some Preliminary Results}\label{s}
Now we recall some results of linear diffusion process.  Let $(X_t^*)_{t \ge 0}$ be a linear diffusion process as $X_t^* = X_0^* + \mu t + \sigma B_t$, where $B_t$ is a standard Brownian motion. Let $\tau_a^* := \inf\{t \ge 0: X_t^* = a\}$ and
\beqnn
l = \frac{\sqrt{2q\sigma^2 + \mu^2}}{\sigma^2}, \quad \delta^+ = l + \frac{\mu}{\sigma^2}, \quad \delta^- = l - \frac{\mu}{\sigma^2}.
\eeqnn
The exponential random variable $e_q$ (rate $q> 0$) is assumed independently of $X^*$. 
For any $x, a \in \mbb{R},$ by 1.0.5 on page 250 and 2.0.1,2.0.2 on page 295 of \cite{BS02}, see also \cite[Lemma 2.1]{CWZ24}, we have
\beqlb\label{star1}
\mbf{P}_x\{X_{e_q}^* \in \d z\} = \frac{q}{l\sigma^2}\exp\left\{\mu(z - x)/\sigma^2 - |z - x|l\right\}\d z,
\eeqlb
\beqlb\label{star2}
\mbf{P}_x\{\tau_a^* < e_q\} = \mbf{E}_x[e^{-q\tau_a^*}] = \exp\left\{\mu(a - x)/\sigma^2 - |a - x|l\right\}
\eeqlb
and
\beqlb\label{star3}
\mbf{P}_x\{\tau_a^* \in \d t\} = \frac{|a - x|}{\sigma \sqrt{2\pi t^3}}\exp\left\{-\frac{(a - x - \mu t)^2}{2 \sigma^2 t}\right\}\d t, \quad t > 0.
\eeqlb
By \cite[(3.10)]{CWZ24}, we have
\beqlb\label{1}
\ar\ar \mbf{P}_x\{X^*_{e_q}\in \d z, \tau_a^* > e_q\} \cr 
\ar\ar\quad = \frac{q}{l\sigma^2}\begin{cases}
	(e^{\delta^-(x - a)} - e^{-\delta^+(x - a)})e^{-\delta^-(z - a)}\d z, & z \ge x \ge a,\\
	(e^{\delta^+(z - a)} - e^{-\delta^-(z - a)})e^{-\delta^+(x - a)}\d z,  & x \ge z \ge a.
\end{cases}
\eeqlb
Moreover, for $a \le x \le b$ we have
\beqlb\label{star4}
  \mbf{P}_x\{\tau_a^* \wedge \tau_b^* < e_q, X^*_{\tau_a^*\wedge\tau_b^*} = a\}\ar=\ar \mbf{E}_x[e^{-q(\tau_a^*\wedge\tau_b^*)}, X^*_{\tau_a^*\wedge\tau_b^*} = a]\cr
\ar=\ar
 e^{\frac{\mu(a - x)}{\sigma^2}}\frac{\sinh((b - x)l)}{\sinh((b - a)l)}
\eeqlb
and
\beqlb\label{star5}
 \mbf{P}_x\{\tau_a^* \wedge \tau_b^* < e_q, X^*_{\tau_a^*\wedge\tau_b^*} = b\} 
\ar =\ar \mbf{E}_x[e^{-q(\tau_a^*\wedge\tau_b^*)}, X^*_{\tau_a^*\wedge\tau_b^*} = b] \cr
\ar=\ar
  e^{\frac{\mu(b - x)}{\sigma^2}}\frac{\sinh((x - a)l)}{\sinh((b - a)l)},
\eeqlb
see 3.05 (a) and 3.05 (b) on page 309 of \cite{BS02}. Moreover, we have the following result.

\begin{lem}\label{L01}
	For any $x \in [a, b]$, we have
	\beqnn
	\mbf{P}_x\{X^*_{e_q} \in \d z, e_q < \tau_a^*\wedge\tau_b^*\} = 
	\begin{cases}
		\frac{2q\sinh((x-a)l)\sinh((b-z)l)}{l\sigma^2 \sinh((b-a)l)}e^{\frac{\mu(z - x)}{\sigma^2}}\d z, & a\le x\le z\le b,\\ 
		\frac{2q \sinh((b-x)l)\sinh((z-a)l)}{l\sigma^2\sinh((b-a)l)}e^{\frac{\mu(z - x)}{\sigma^2}}\d z, & a\le z< x\le b. 
	\end{cases}
	\eeqnn
\end{lem}

\proof
For the case of $a \le x \le z \le b,$ by \eqref{star1}, \eqref{star4}, \eqref{star5} and the strong Markov property, we have
\beqlb\label{eq0305b}
\ar\ar\mbf{P}_x\{X^*_{e_q} \in \d z, e_q < \tau_a^*\wedge\tau_b^*\} \cr
\ar\ar\quad= \mbf{P}_x\{X^*_{e_q} \in \d z\} - \mbf{P}_x\{X^*_{e_q} \in \d z, e_q \ge \tau_a^*\wedge\tau_b^*\}\cr
\ar\ar\quad= \mbf{P}_x\{X^*_{e_q} \in \d z\} - \mbf{P}_x\{X^*_{e_q} \in \d z, X^*_{\tau_a^*\wedge\tau_b^*} = a, e_q \ge \tau_a^*\wedge\tau_b^*\}\cr
\ar\ar\qquad - \mbf{P}_x\{X^*_{e_q} \in \d z, X^*_{\tau_a^*\wedge\tau_b^*} = b, e_q \ge \tau_a^*\wedge\tau_b^*\}\cr
\ar\ar\quad= \mbf{P}_x\{X^*_{e_q} \in \d z\} - \mbf{P}_x\{X_{\tau_a^*\wedge\tau_b^*} = a, e_q \ge \tau_a^*\wedge\tau_b^*\}\mbf{P}_a\{X^*_{e_q} \in \d z\}\cr
\ar\ar\qquad -  \mbf{P}_x\{X_{\tau_a^*\wedge\tau_b^*} = b, e_q \ge \tau_a^*\wedge\tau_b^*\}\mbf{P}_b\{X^*_{e_q} \in \d z\}\cr
\ar\ar\quad= \frac{q}{l\sigma^2}\Bigg[e^{-\delta^-(z - x)} -  e^{\frac{\mu(a - x)}{\sigma^2}-\delta^-(z - a)}\frac{\sinh((b - x)l)}{\sinh((b - a)l)}\cr
\ar\ar\qquad\qquad\quad - e^{\frac{\mu(b - x)}{\sigma^2}+\delta^+(z - b)}\frac{\sinh((x - a)l)}{\sinh((b - a)l)}\Bigg]\d z.
\eeqlb
Notice that
\beqnn
\ar\ar e^{-\delta^-(z - x)} \sinh\left((b - a)l\right)-  e^{\frac{\mu(a - x)}{\sigma^2}-\delta^-(z - a)}\sinh\left((b - x)l\right)  - e^{\frac{\mu(b - x)}{\sigma^2}+\delta^+(z - b)}\sinh\left((x - a)l\right)\cr
\ar\ar\quad= \frac{1}{2}e^{\frac{\mu(z - x)}{\sigma^2}}\Big[e^{l(b-a-z+x)}- e^{l(b-x-z+a)} - e^{l(x-a+z-b)} + e^{l(z-b-x+a)}\Big]\cr
\ar\ar\quad = 2e^{\frac{\mu(z - x)}{\sigma^2}}\sinh((x - a)l)\sinh((b-z)l).
\eeqnn
Then by the above and \eqref{eq0305b}, for $a\le x\le z \le b$, we have
\beqlb\label{eq0307b}
\mbf{P}_x\{X^*_{e_q} \in \d z, e_q < \tau_a^*\wedge\tau_b^*\}  = \frac{2q\sinh((x-a)l)\sinh((b-z)l)}{l\sigma^2 \sinh((b-a)l)}e^{\frac{\mu(z - x)}{\sigma^2}}\d z.  
\eeqlb

Similarly, for $a\le z < x\le b$, we have 
\beqlb\label{eq0307c}
\ar\ar \mbf{P}_{x}\{X^*_{e_q} \in \d z, e_q < \tau_a^*\wedge\tau_b^*\} \cr
\ar\ar\quad = \mbf{P}_{x}\{X^*_{e_q} \in \d z\} - \mbf{P}_{x}\{ X_{\tau_a^*\wedge\tau_b^*} = a, e_q \ge \tau_a^*\wedge\tau_b^*\}\mbf{P}_a\{X^*_{e_q} \in \d z\}\cr
\ar\ar\qquad - \mbf{P}_{x}\{ X_{\tau_a^*\wedge\tau_b^*} = b, e_q \ge \tau_a^*\wedge\tau_b^*\}\mbf{P}_b\{X^*_{e_q} \in \d z\}\cr
\ar\ar\quad = \frac{q}{l\sigma^2}\Bigg[e^{\delta^+(z - x)} -  e^{\frac{\mu(a - x)}{\sigma^2}-\delta^-(z - a)}\frac{\sinh(l(b - x))}{\sinh(l(b - a))}\cr
\ar\ar\qquad\qquad\quad - e^{\frac{\mu(b - x)}{\sigma^2}+\delta^+(z - b)}\frac{\sinh(l(x - a))}{\sinh(l(b - a))}\Bigg]\d z\cr
\ar\ar\quad= \frac{2q \sinh((b-x)l)\sinh((z-a)l)}{l\sigma^2\sinh((b-a)l)}e^{\frac{\mu(z - x)}{\sigma^2}}\d z.
\eeqlb
The result follows from \eqref{eq0307b} and \eqref{eq0307c}.
\qed

Now we specify the constants $b_i^+, c^+_{i}, b_i^-$ and $c^-_{i}$ in \eqref{g_q^+} and \eqref{g_q^-}. Let 
\beqlb\label{b^+}
\begin{cases}
	b_i^+ = b_{i-1}^+\left[(1-c^+_{i-1})e^{\delta_{i-1}^-(a_{i} - a_{i-1})} + c^+_{i-1}e^{-\delta_{i-1}^+(a_{i} - a_{i-1})}\right], \  i = 2, \cdots, n,\\
	b_1^+ = 1
\end{cases}
\eeqlb
and
\beqlb\label{c^+_i}
\begin{cases}
	c^+_{i} = \frac{\delta_i^- - \delta_{i-1}^-}{\delta_i^+ + \delta_i^-} + \frac{\delta_{i-1}^+ + \delta_{i-1}^-}{\delta_i^+ + \delta_i^-} \frac{c^+_{i-1}}{(1-c^+_{i-1})e^{(\delta_{i-1}^- + \delta_{i - 1}^+)(a_{i} - a_{i-1})} + c^+_{i-1}}, \   i = 2, \cdots, n,\\
	c^+_{1} = \frac{\delta_1^- - \delta_0^-}{\delta_1^+ + \delta_1^-}.
\end{cases}
\eeqlb
Let
\beqlb\label{b^-}
\begin{cases}
	b_i^- = b_{i+1}^- \left[c^-_{i+1}e^{-\delta_i^-(a_{i+1} - a_i)} + (1- c^-_{i+1})e^{\delta_i^+(a_{i+1} - a_i)}\right], \   i = 1, \cdots, n-1,\\
	b_n^- = 1
\end{cases}
\eeqlb
and
\beqlb\label{c_i^-}
\begin{cases}
	c^-_{i} =  \frac{\delta_{i-1}^+ - \delta_i^+}{\delta_{i-1}^+ + \delta_{i-1}^-} + \frac{\delta_i^- + \delta_i^+}{\delta_{i-1}^+ + \delta_{i-1}^-}\frac{ c^-_{i+1}}{c^-_{i+1} + (1- c^-_{i+1})e^{(\delta_i^+ + \delta_i^-)(a_{i+1} - a_i)}}, \   i = 1, \cdots, n -1,\\
	c^-_{n} = \frac{\delta_{n-1}^+ - \delta_n^+}{\delta_{n-1}^+ + \delta_{n-1}^-}.
\end{cases}
\eeqlb
Then we have the following result.
\begin{lem}
	For constants $b_i^+, c^+_{i}, b_i^-$ and $c^-_{i}$ given by \eqref{b^+}-\eqref{c_i^-}, we have $g_q^{\pm} \in C^1(\mbb{R})\cap C^2(\mbb{R}\setminus \cup_{i = 1}^n\{a_i\})$.
\end{lem}
 \proof
 By \eqref{b^+} and \eqref{b^-}, we have $g_q^{+}(a_1-) = 1 = g_q^{+}(a_1+) = b_1^+$
and for $i = 2, \cdots, n,$
\beqnn
g_q^{+}(a_i-) \ar=\ar b_{i-1}^+ \left[(1-c^+_{i-1})e^{\delta_{i-1}^-(a_{i} - a_{i-1})} + c^+_{i-1}e^{-\delta_{i-1}^+(a_{i} - a_{i-1})}\right]\cr
\ar=\ar g_q^{+}(a_i+) = b_i^+.
\eeqnn
Moreover, one sees that $g_q^{-}(a_n+) = 1 = g_q^{-}(a_n-) = b_n^-$ and for $i = 1, \cdots, n - 1,$ 
\beqnn
g_q^{-}(a_i+) \ar=\ar b_{i+1}^-\left[c^-_{i+1}e^{\delta_i^-(a_i - a_{i+1})} + (1- c^-_{i+1})e^{-\delta_i^+(a_i - a_{i+1})}\right]\cr
\ar=\ar g_q^{-}(a_i-) = b_i^-,
\eeqnn
Then it follows that $g_q^{\pm}(a_i-) = g_q^{\pm}(a_i+)$ for $i = 1, \cdots, n$. 

On the other hand, by \eqref{b^+}-\eqref{c_i^-}, we have
\beqnn
g_q^{+'}(a_1-) = \delta_0^- = g_q^{+'}(a_1+) = b_1^+\left[(1 - c^{+}_{1})\delta_1^- - c^{+}_{1}\delta_1^+\right],
\eeqnn 
and for $i = 2, \cdots, n,$  
\beqnn
g_q^{+'}(a_i-) \ar=\ar b_{i-1}^+\left[(1-c^+_{i-1})\delta_{i-1}^-e^{\delta_{i-1}^-(a_{i} - a_{i-1})} - c^+_{i-1}\delta_{i-1}^+e^{-\delta_{i-1}^+(a_{i} - a_{i-1})}\right]\cr
\ar=\ar g_q^{+'}(a_i+) = b_i^+\left[(1-c^+_{i})\delta_i^- - c^+_{i}\delta_i^+\right].
\eeqnn 
Moreover, one sees that
\beqnn
g_q^{-'}(a_n+) = -\delta_n^+ = g_q^{-'}(a_n-) = b_n^-\left[c^-_{n}\delta_{n-1}^- - (1-c^-_{n})\delta_{n-1}^+\right],
\eeqnn
and for $i = 1, \cdots, n -1,$
\beqnn
g_q^{-'}(a_i+) \ar=\ar b_{i+1}^-\left[c^-_{i+1}\delta_{i}^-e^{\delta_i^-(a_i - a_{i+1})} - (1- c^-_{i+1})\delta_{i}^+e^{-\delta_i^+(a_i - a_{i+1})}\right]\cr
\ar=\ar g_q^{-'}(a_i-) = b_i^-\left[c^-_{i}\delta_{i-1}^- - (1 - c^-_{i})\delta_{i-1}^+\right].
\eeqnn
Then $g_q^{\pm '}(a_i-) = g_q^{\pm '}(a_i+)$ for $i = 1, \cdots, n$. The result follows.
\qed

\section{Proof of Theorem \ref{tpm} }\label{S3}

In this section, we give the proofs of Theorem \ref{tpm}. Recall that $\tau_{a_i} = \inf\{t \ge 0: X_t = a_i\}$ for $i = 0, 1, \cdots, n.$ We use the simplified notation $\tau_i$ for $\tau_{a_i}$ for convenience. let $X^i:=(X_t^i)_{t \ge 0}$ be defined as
\beqnn
 X_t^i = \mu_i t + \sigma_i B_t, \qquad t\geq 0, \, i = 0, \cdots, n.
\eeqnn
For initial value $x < a_1$ (resp. $x > a_n$), $X_t = X_t^0$ for any $t \in [0, \tau_1)$ (resp. $X_t = X_t^n$ for any $t \in [0, \tau_n)$) almost surely. Moreover, the stopping time $\tau_1$ (resp. $\tau_n$) coincides with the first passage time of $X^0$ (resp. $X^n$) to $a_1$ (resp. $a_n$).  Similarly, for $x \in (a_i, a_{i+1})$ with $i = 1, \cdots, n-1$, $\mbf{P}_x\{X_t = X_t^i\ \text{for\ any}\ t \in [0, \tau_i \wedge \tau_{i+1})) = 1$, and the stopping time $\tau_i \wedge \tau_{i+1}$ corresponds to the first exit time of $X^i$ from the interval $(a_i, a_{i+1}).$ The properties of $X^i, i = 0, \cdots, n$ can then be used to derive those of $X$.

\subsection{The case  $z \le a_1$}

\begin{lem}\label{p0805c}
	For any $z \le a_1$, we have
	\beqnn
	\mbf{P}_{a_1}\{X_{e_q}\in \d z\}=  \frac{q e^{\delta_0^+(z - a_1)}}{l_0\sigma_0^2(1-c^-_{1})}\d z.
	\eeqnn
\end{lem}
\proof 
For any $r \in (z, a_1),$ one obtains by the Markov property that
\beqnn
\mbf{P}_{a_1}\{X_{e_q} \in \d z\} \ar=\ar\mbf{P}_{a_1}\{X_{e_q} \in \d z, \tau_r \le e_q\}\cr
\ar=\ar \mbf{P}_{a_1}\{\tau_r \le e_q\}\mbf{P}_r\{X_{e_q} \in \d z\}\cr
\ar=\ar \mbf{E}_{a_1}[e^{-q\tau_r}]\left[\mbf{P}_r\{X_{e_q} \in \d z, \tau_1 > e_q\} + \mbf{P}_r\{X_{e_q} \in \d z, \tau_1 \le e_q\}\right]\cr
\ar=\ar \mbf{E}_{a_1}[e^{-q\tau_r}]\left[\mbf{P}_r\{X^1_{e_q} \in \d z, \tau_1 > e_q\} + \mbf{P}_r\{\tau_1 \le e_q\}\mbf{P}_{a_1}\{X_{e_q} \in \d z\}\right]\cr
\ar=\ar \mbf{E}_{a_1}[e^{-q\tau_r}]\left[\mbf{P}_r\{X^1_{e_q} \in \d z, \tau_1 > e_q\} + \mbf{E}_r[e^{-q\tau_1}]\mbf{P}_{a_1}\{X_{e_q} \in \d z\}\right].
\eeqnn
Then
\beqlb\label{eq2}
\mbf{P}_{a_1}\{X_{e_q} \in \d z\} = \frac{\mbf{E}_{a_1}[e^{-q\tau_r}]\mbf{P}_r\{X^1_{e_q} \in \d z, \tau_1 > e_q\}}{1 - \mbf{E}_{a_1}[e^{-q\tau_r}]\mbf{E}_r[e^{-q\tau_1}]}.
\eeqlb
Notice that $\tau_1 = \inf\{t \ge 0: X_t^1 = a_1\}$ when $X_0 = X_0^1 = r < a_1.$ It follows by \cite[(3.14)]{CWZ24} that
\beqnn
\mbf{P}_r\{X^1_{e_q} \in \d z, \tau_1 > e_q\} = \frac{q}{l_0\sigma_0^2}\left(e^{-\delta_0^+(r - a_1)} - e^{\delta_0^-(r - a_1)}\right)e^{\delta_0^+(z - a_1)}\d z.
\eeqnn
Moreover, by Corollary \ref{exptau} we have
\beqnn
\mbf{E}_{a_1}[e^{-q\tau_r}] = \frac{g_q^-(a_1)}{g_q^-(r)} = \frac{1}{c^-_{1}e^{\delta_{0}^-(r-a_1)}+(1-c^-_{1})e^{-\delta_{0}^+(r-a_1)}}
\eeqnn
and
\beqnn
\mbf{E}_r[e^{-q\tau_1}] = \frac{g_q^+(r)}{g_q^+(a_1)} = e^{\delta_0^-(r - a_1)}.
\eeqnn
Then by the above and \eqref{eq2}, one obtains that
\beqnn
\mbf{P}_{a_1}\{X_{e_q} \in \d z\} \ar=\ar \frac{\mbf{P}_r\{X^1_{e_q} \in \d z, \tau_1 > e_q\}}{\frac{1}{\mbf{E}_{a_1}[e^{-q\tau_r}]} - \mbf{E}_r[e^{-q\tau_1}]}\cr
\ar=\ar \frac{\frac{q}{l_0\sigma_0^2}(e^{-\delta_0^+(r - a_1)} - e^{\delta_0^-(r - a_1)})e^{\delta_0^+(z - a_1)}\d z}{c^-_{1}e^{\delta_{0}^-(r-a_1)}+(1-c^-_{1})e^{-\delta_{0}^+(r-a_1)} -  e^{\delta_0^-(r - a_1)}}\cr
\ar=\ar \frac{q e^{\delta_0^+(z - a_1)}}{l_0\sigma_0^2(1-c^-_{1})}\d z.
\eeqnn
The proof is completed.
\qed

\begin{prop}\label{p0717b}
	For any $x, z \le a_1$, we have
	\beqnn
	\mbf{P}_{x}\{X_{e_q}\in \d z\}= \begin{cases}
	 \frac{q e^{-\delta_0^+(a_1 - z)}/\mbf{E}_{a_1} [e^{-q\tau_x}]}{l_0\sigma_0^2 (1-c^-_1)} \d z, &   a_1 > x\ge z,\cr
	e^{2\mu_0(z - x)/\sigma_0^2}\frac{q e^{-\delta_0^+(a_1 - x)}/\mbf{E}_{a_1}[e^{-q\tau_z}]}{l_0\sigma_0^2(1-c^-_{1})} \d z, &   a_1 > z\ge x.
\end{cases}
	\eeqnn
\end{prop}

\proof	
For any $x, z \le a_1,$ by the Markov property we have
\beqlb\label{eq0311}
\mbf{P}_{x}\{X_{e_q}\in \d z\}\ar=\ar \mbf{P}_{x}\{X_{e_q}\in \d z, e_q \le \tau_1\} + \mbf{P}_x\{X_{e_q}\in \d z, e_q > \tau_1\}\cr
\ar=\ar \mbf{P}_{x}\{X_{e_q}^1\in \d z, e_q \le \tau_1 \} + \mbf{P}_x\{\tau_1 < e_q\}\mbf{P}_{a_1}\{X_{e_q}\in \d z\}\cr
\ar=\ar  \mbf{P}_{x}\{X_{e_q}^1\in \d z, e_q \le \tau_1 \} + \mbf{E}_x[e^{-q\tau_1}]\mbf{P}_{a_1}\{X_{e_q}\in \d z\}.
\eeqlb
Again, notice that $\tau_1 = \inf\{t \ge 0: X_t^1 = a_1\}$ with $X_0 = X_0^1 = x < a_1.$ Then by \eqref{1}, one obtains that
\beqlb\label{eq0311a}
\ar\ar\mbf{P}_{x}\{X_{e_q}^1\in \d z, e_q \le \tau_1 \} \cr
\ar\ar\quad = \begin{cases}
	\frac{q}{l_0\sigma_0^2}(e^{\delta_0^+( a_1-x)} - e^{-\delta_0^-( a_1-x)})e^{\delta_0^+(z - a_1)}\d z, & a_1\ge x\ge z,\cr
	\frac{q}{l_0\sigma_0^2}(e^{\delta_0^-( a_1-z)} - e^{-\delta_0^+( a_1-z)})e^{\delta_0^-(x - a_1)}\d z, & a_1\ge z\ge x.
\end{cases}
\eeqlb
On the other hand, by Corollary \ref{exptau}, we have
\beqlb\label{eq0311b}
\mbf{E}_x[e^{-q\tau_1}] = \frac{g_q^-(x)}{g_q^-(a_1)} = e^{\delta_0^-(x - a_1)} 
\eeqlb
and
\beqlb\label{l0809a}
\mbf{E}_{a_1} [e^{-q\tau_x}] = \frac{1}{(1 - c_1^-)e^{\delta_0^+(a_1 - x)} + c_1^- e^{-\delta_0^-(a_1 - x)}}.
\eeqlb
It follows from \eqref{eq0311}-\eqref{eq0311b} and Lemma \ref{p0805c} that 
\beqnn
 \mbf{P}_{x}\{X_{e_q}\in \d z\}  
\ar=\ar \begin{cases}
	\frac{q e^{\delta_0^+(z - a_1)}}{l_0\sigma_0^2} \big[e^{\delta_0^+( a_1-x)} - e^{-\delta_0^-( a_1-x)}+\frac{e^{\delta_0^-(x - a_1)}}{1-c^-_{1}}\big]  \d z, & a_1 \ge x\ge z,\cr
	\frac{q e^{\delta_0^-(x - a_1)}}{l_0\sigma_0^2} \big[e^{\delta_0^-( a_1-z)} - e^{-\delta_0^+( a_1-z)}+\frac{e^{\delta_0^+(z - a_1)}}{1-c^-_{1}}\big] \d z, &  a_1 \ge z\ge x
\end{cases}\cr
\ar=\ar \begin{cases}
	e^{\mu_0(z - x)/\sigma_0^2}\frac{q e^{-l_0(a_1 - z)}}{l_0\sigma_0^2} \left[e^{l_0( a_1-x)} +\frac{c^-_{1}e^{-l_0(a_1 - x)}}{1-c^-_{1}}\right]  \d z, & a_1 \ge x\ge z,\cr
	e^{\mu_0(z - x)/\sigma_0^2}\frac{q e^{-l_0(a_1 - x)}}{l_0\sigma_0^2} \left[e^{l_0( a_1-z)} +\frac{c^-_{1}e^{-l_0(a_1 - z)}}{1-c^-_{1}}\right] \d z, & a_1 \ge z\ge x.
\end{cases}
\eeqnn
Then the proof is completed by the above and \eqref{l0809a}. 
\qed

\begin{prop}\label{p0805d}
	For any $x \ge a_1$ and $z \le a_1$, we have
	\beqnn
	\mbf{P}_x\{X_{e_q}\in \d z\} = \frac{q e^{\delta_0^+(z - a_1)}\mbf{E}_x[e^{-q\tau_1}]}{l_0\sigma_0^2(1-c^-_{1})}\d z.
	\eeqnn
\end{prop}
\proof
For the case of $x \ge a_1 \ge z$, by \eqref{g_q^-} and the strong Markov property, one sees that
\beqnn
\mbf{P}_x\{X_{e_q}\in \d z\} \ar=\ar \mbf{P}_x\{X_{e_q}\in \d z, \tau_1 < e_q\} = \mbf{E}_x[e^{-q\tau_1}]\mbf{P}_{a_1}\{X_{e_q}\in \d z\}.
\eeqnn
The result follows from Proposition \ref{p0717b}.
\qed

\subsection{The case $z \in (a_i, a_{i+1}]$ for $i = 1, \cdots, n-1$}

\begin{lem}\label{l0801a}
Let $X$ be the solution to \eqref{Xn}.
For any $z \in (a_i, a_{i+1}]$ and $i = 1, \cdots, n-1,$ we have
\beqnn
\mbf{P}_{a_i}\{X_{e_q} \in \d z, e_q < \tau_{i+1} \}  =  \frac{2q\sinh((a_{i+1}-z)l_i)}{l_i\sigma_i^2}e^{-\mu_i(a_{i+1} - z)/\sigma_i^2} \mbf{E}_{a_i}[e^{-q\tau_{i+1}}] \d z.
\eeqnn
\end{lem}
\proof
For any $r \in (a_i, z)$, one sees that $\tau_i\wedge \tau_{i+1} = \inf\{t \ge 0: X_t^i \notin (a_i, a_{i+1})\}$ with $X_0 = X_0^i = r \in [a_i, a_{i+1}].$ Then by the Markov property, we have
\beqnn
\ar\ar \mbf{P}_r\{X_{e_q} \in \d z, \tau_{i} < e_q < \tau_{i + 1}\}\cr
\ar\ar\quad= \mbf{P}_r\{X_{e_q} \in \d z, \tau_i\wedge \tau_{i+1} < e_q, X_{ \tau_i\wedge \tau_{i+1}} = a_i, e_q < \tau_{i + 1}\}\cr
\ar\ar\quad=  \mbf{P}_r\{\tau_i\wedge \tau_{i+1} < e_q, X_{\tau_i\wedge \tau_{i+1}} = a_i\}\mbf{P}_{a_i} \{X_{e_q} \in \d z, e_q < \tau_{i + 1}\}\cr
\ar\ar\quad=  \mbf{P}_r\{\tau_i \wedge \tau_{i+1}  < e_q, X^i_{\tau_i\wedge \tau_{i+1}} = a_i\}\mbf{P}_{a_i} \{X_{e_q} \in \d z, e_q < \tau_{i + 1}\}
\eeqnn
and
\beqnn
\ar\ar\mbf{P}_{a_i}\{X_{e_q} \in \d z, \tau_r < e_q < \tau_{i+1}\}\cr
 \ar\ar\quad=  \mbf{P}_{a_i}\{\tau_r\wedge \tau_{i+1} < e_q, X_{\tau_r\wedge \tau_{i+1}} = r\}\mbf{P}_{r} \{X_{e_q} \in \d z, e_q < \tau_{i + 1}\}\cr
\ar\ar\quad=  \mbf{P}_{a_i}\{\tau_r < e_q\} \mbf{P}_{r} \{X_{e_q} \in \d z, e_q < \tau_{i + 1}\},
\eeqnn
which implies  
\beqnn
\ar\ar \mbf{P}_{a_i}\{X_{e_q} \in \d z, e_q < \tau_{i+1}\}\cr
\ar\ar \quad =  \mbf{P}_{a_i}\{X_{e_q} \in \d z, \tau_r < e_q < \tau_{i+1}\}\cr
\ar \ar\quad = \mbf{P}_{a_i}\{\tau_r < e_q\}\mbf{P}_r\{X_{e_q} \in \d z, e_q < \tau_{i+1}\}\cr
\ar \ar\quad = \mbf{E}_{a_i}[e^{-q\tau_r}] \Big[\mbf{P}_r\{X_{e_q} \in \d z, e_q < \tau_{i+1}\wedge\tau_i \} + \mbf{P}_r\{X_{e_q} \in \d z, \tau_{i} < e_q < \tau_{i + 1}\}\Big]\cr
\ar\ar\quad =\mbf{E}_{a_i}[e^{-q\tau_r}]\mbf{P}_r\{X_{e_q}^i \in \d z, e_q < \tau_{i+1} \wedge\tau_i\}\cr
\ar\ar \qquad + \mbf{E}_{a_i}[e^{-q\tau_r}] \mbf{P}_r\{\tau_i \wedge \tau_{i+1} < e_q, X^i_{\tau_i\wedge \tau_{i+1}} = a_i\}\mbf{P}_{a_i}\{X_{e_q} \in \d z, e_q < \tau_{i+1}\}. 
\eeqnn
It follows that
\beqlb\label{eq337a}
\mbf{P}_{a_i}\{X_{e_q} \in \d z, e_q < \tau_{i+1}\} \ar= \ar \frac{\mbf{E}_{a_i}[e^{-q\tau_r}]\mbf{P}_r\{X_{e_q}^i \in \d z, e_q < \tau_{i+1}\wedge\tau_i\}}{1 -  \mbf{E}_{a_i}[e^{-q\tau_r}] \mbf{P}_r\{\tau_i\wedge \tau_{i+1} < e_q, X^i_{\tau_i\wedge \tau_{i+1}} = a_i\}}\cr
\ar= \ar \frac{\mbf{P}_r\{X_{e_q}^i \in \d z, e_q <  \tau_{i+1}\wedge\tau_i\}}{\frac{1}{\mbf{E}_{a_i}[e^{-q\tau_r}]} -  \mbf{P}_r\{\tau_i\wedge \tau_{i+1} < e_q, X^i_{\tau_i\wedge \tau_{i+1}} = a_i\}}.
\eeqlb
By Lemma \ref{L01} we have
\beqlb\label{eq337}
\ar\ar \mbf{P}_r\{X_{e_q}^i \in \d z, e_q < \tau_{i+1}\wedge\tau_i\} \cr
\ar\ar\quad = \frac{2q\sinh((r-a_i)l_i)\sinh((a_{i+1}-z)l_i)}{l_i\sigma_i^2 \sinh((a_{i+1}-a_i)l_i)}e^{\mu_i(z - r)/\sigma_i^2}\d z.
\eeqlb
For $i = 1, \cdots, n-1,$ by Corollary \ref{exptau} and Proposition \eqref{star4}, one obtains
\beqnn
\mbf{E}_{a_i}[e^{-q\tau_r}] \ar=\ar \frac{g_q^+(a_i)}{g_q^+(r)} = \frac{1}{(1 - c_i^+)e^{\delta_{i}^-(r - a_{i})} + c_i^+e^{-\delta_{i}^+(r - a_{i})}} 
\eeqnn
and
\beqnn
\mbf{P}_r\{\tau_i \wedge \tau_{i+1} < e_q, X^i_{\tau_i^i\wedge\tau_{i+1}^i} = a_i\}
= e^{\frac{\mu_i(a_i - r)}{\sigma_i^2}}\frac{\sinh((a_{i+1} - r)l_i)}{\sinh((a_{i+1} - a_i)l_i)}.
\eeqnn
It follows that
\beqnn
\ar\ar 2e^{\frac{\mu_i(r - a_i)}{\sigma_i^2}}\sinh((a_{i+1} - a_i)l_i)\left[\frac{1}{\mbf{E}_{a_i}[e^{-q\tau_r}]}-\mbf{P}_r\{\tau_i^i \wedge \tau_{i+1}^i < e_q, X^i_{\tau_i^i\wedge\tau_{i+1}^i} = a_i\}\right]\cr
\ar\ar\quad= 2\sinh((a_{i+1} - a_i)l_i)\left[(1 - c_i^+)e^{l_i(r - a_i)} + c_i^+ e^{-l_i(r - a_i)}\right] - 2\sinh((a_{i+1} - r)l_i) \cr
\ar\ar\quad=(1 - c_i^+)e^{l_i(a_{i+1} - a_i)}\left[e^{l_i(r - a_i)} - e^{-l_i(r - a_i)}\right] + c_i^+e^{-l_i(a_{i + 1} - a_i)}\left[e^{l_i(r - a_i)} - e^{-l_i(r - a_i)}\right]\cr
\ar\ar\quad = 2\sinh(l_i(r - a_i)) \left[(1 - c_i^+)e^{l_i(a_{i+1} - a_i)}   + c_i^+e^{-l_i(a_{i + 1} - a_i)}\right].
\eeqnn
By the above, \eqref{eq337a} and \eqref{eq337}, we have
\beqnn
\mbf{P}_{a_i}\{X_{e_q} \in \d z, \tau_{i+1} > e_q\} 
\ar=\ar \frac{2q\sinh((a_{i+1}-z)l_i)}{l_i\sigma_i^2} \frac{e^{\mu_i(z - a_i)/\sigma_i^2}\d z}{(1 - c_i^+)e^{l_i(a_{i+1} - a_i)}   + c_i^+e^{-l_i(a_{i + 1} - a_i)}}\cr
\ar=\ar  \frac{2q\sinh((a_{i+1}-z)l_i)}{l_i\sigma_i^2}e^{-\mu_i(a_{i+1} - z)/\sigma_i^2} \mbf{E}_{a_i}[e^{-q\tau_{i+1}}]\d z.
\eeqnn
The result follows.
\qed

\begin{lem}\label{l0801b}
For any $z \in (a_i, a_{i+1}]$ and $i = 1, \cdots, n - 1,$ we have 
\beqnn
\mbf{P}_{a_{i+1}}\{X_{e_q} \in \d z, e_q < \tau_{i}\} = \frac{2q \sinh((z-a_{i})l_i)}{l_i\sigma_i^2}e^{\mu_i(z - a_{i})/\sigma_i^2}\mbf{E}_{a_{i+1}}[e^{-q\tau_i}] \d z.
\eeqnn
\end{lem}
\proof
The proof is similar to the one of Lemma \ref{l0801a}. By the strong Markov property, for $r \in (z, a_{i+1}),$ we have
\beqnn
\ar\ar \mbf{P}_{a_{i+1}}\{X_{e_q} \in \d z, e_q < \tau_{i}\} \cr
\ar\ar\quad = \mbf{P}_{a_{i+1}}\{X_{e_q} \in \d z,\tau_{r}< e_q < \tau_{i}\} = \mbf{P}_{a_{i+1}}\{\tau_r <e_q\}\mbf{P}_r\{X_{e_q}\in\d z,  e_q < \tau_i\}\cr
\ar\ar\quad = \mbf{E}_{a_{i+1}}[e^{-q\tau_r}] \mbf{P}_r\{X_{e_q}\in\d z,  e_q < \tau_i\}\cr
\ar\ar\quad = \mbf{E}_{a_{i+1}}[e^{-q\tau_r}] \Big[\mbf{P}_r\{X_{e_q} \in \d z, \tau_{i+1}\wedge\tau_i > e_q\} + \mbf{P}_r\{X_{e_q} \in \d z, \tau_{i} > e_q > \tau_{i+1}\}\Big]\cr
\ar\ar\quad = \mbf{E}_{a_{i+1}}[e^{-q\tau_r}]\Big[\mbf{P}_r\{X_{e_q}^{i} \in \d z, \tau_{i+1}\wedge\tau_i > e_q\}\cr
\ar\ar\qquad  + \mbf{P}_r\{\tau_i\wedge\tau_{i + 1} < e_q, X^i_{ \tau_i\wedge\tau_{i + 1}} = a_{i + 1}\}\mbf{P}_{a_{i+1}}\{X_{e_q} \in \d z, \tau_{i} > e_q\}\Big],
\eeqnn
which implies that
\beqlb\label{eq37c}
\mbf{P}_{a_{i+1}}\{X_{e_q} \in \d z, e_q<\tau_{i}  \} \ar=\ar \frac{\mbf{E}_{a_{i+1}}[e^{-q\tau_r}]\mbf{P}_r\{X_{e_q}^{i} \in \d z, \tau_{i+1}\wedge\tau_i > e_q\}}{1 -  \mbf{E}_{a_{i+1}}[e^{-q\tau_r}]\mbf{P}_r\{\tau_i\wedge\tau_{i + 1} < e_q, X^i_{ \tau_i\wedge\tau_{i + 1}} = a_{i + 1}\}}\cr
\ar=\ar \frac{\mbf{P}_r\{X_{e_q}^{i} \in \d z, \tau_{i+1}\wedge\tau_i > e_q\}}{\frac{1}{\mbf{E}_{a_{i+1}}[e^{-q\tau_r}]} - \mbf{P}_r\{\tau_i\wedge\tau_{i + 1} < e_q, X^i_{ \tau_i\wedge\tau_{i + 1}} = a_{i + 1}\}}.
\eeqlb
By Lemma \ref{L01}, one sees that
\beqlb\label{eq37d}
\ar\ar \mbf{P}_r\{X_{e_q}^i \in \d z, \tau_{i+1}\wedge\tau_i > e_q\} \cr
\ar\ar\quad =\frac{2q \sinh((a_{i+1}-r)l_i)\sinh((z-a_{i})l_i)}{l_i\sigma_i^2\sinh((a_{i+1}-a_{i})l_i)}e^{\mu_i(z - r)/\sigma_i^2}\d z.
\eeqlb
For $i = 1, \cdots, n-1,$  by Corollary \ref{exptau} and Theorem \eqref{star5} we have 
\beqnn
\mbf{E}_{a_{i+1}}[e^{-q\tau_r}] = \frac{g_q^-(a_{i+1})}{g_q^-(r)} \ar=\ar \frac{1}{c_{i+1}^-e^{\delta_{i}^-(r - a_{i+1})} + (1 - c_{i+1}^-)e^{-\delta_{i}^+(r - a_{i+1})}} 
\eeqnn
and
\beqnn
\mbf{P}_r\{\tau_i\wedge\tau_{i + 1} < e_q, X^i_{ \tau_i\wedge\tau_{i + 1}} = a_{i + 1}\} = \frac{\sinh((r - a_i)l_i)}{\sinh((a_{i + 1} - a_i)l_i)}e^{\frac{\mu_i(a_{i + 1} - r)}{\sigma_i^2}}.
\eeqnn
It follows that
\beqnn
\ar\ar 2e^{-\frac{\mu_i(a_{i+1} - r)}{\sigma_i^2}} \sinh((a_{i + 1} - a_i)l_i)\left[ \frac{1}{\mbf{E}_{a_{i+1}}[e^{-q\tau_r}]} - \mbf{P}_r\{\tau_i^i\wedge\tau^i_{i + 1} < e_q, X^i_{ \tau_i\wedge\tau_{i + 1}} = a_{i + 1}\}\right]\cr
\ar\ar\quad = 2\sinh((a_{i+1} - a_i)l_i) [c_{i+1}^-e^{l_{i}(r - a_{i+1})} + (1 - c_{i+1}^-)e^{-l_{i}(r - a_{i+1})}] - 2\sinh((r - a_i)l_i) \cr
\ar\ar\quad = -(1 - c_{i+1}^-)e^{l_{i}(r - a_{i})} + (1 - c_{i+1}^-)e^{l_{i}(2a_{i+1} - a_i - r)} - c_{i+1}^-e^{-l_{i}(2a_{i+1} - a_i - r)} + c_{i+1}^-e^{-l_{i}(r - a_{i})}\cr
\ar\ar\quad= 2\sinh(l_i(a_{i+1} - r))\left[(1 - c_{i+1}^-)e^{l_{i}(a_{i+1} - a_i)} + c_{i+1}^-e^{-l_{i}(a_{i+1} - a_{i})}\right].
\eeqnn
Combining the above, \eqref{eq37d} and \eqref{eq37c}, one obtains that
\beqnn
\ar\ar \mbf{P}_{a_{i+1}}\{X_{e_q} \in \d z, e_q<\tau_{i}  \} \cr
\ar\ar\quad = \frac{2q \sinh((z-a_{i})l_i)}{l_i\sigma_i^2}e^{\mu_i(z - a_{i+1})/\sigma_i^2} \frac{\d z}{(1 - c_{i+1}^-)e^{l_{i}(a_{i+1} - a_i)} + c_{i+1}^-e^{-l_{i}(a_{i+1} - a_{i})}}\cr
\ar\ar\quad = \frac{2q \sinh((z-a_{i})l_i)}{l_i\sigma_i^2}e^{\mu_i(z - a_{i})/\sigma_i^2} \frac{\d z}{(1 - c_{i+1}^-)e^{\delta_{i}^+(a_{i+1} - a_i)} + c_{i+1}^-e^{-\delta^-_{i}(a_{i+1} - a_{i})}}\cr
\ar\ar\quad = \frac{2q \sinh((z-a_{i})l_i)}{l_i\sigma_i^2}e^{\mu_i(z - a_{i})/\sigma_i^2}\mbf{E}_{a_{i+1}}[e^{-q\tau_i}]\d z.
\eeqnn
The result follows.
\qed

\begin{lem}\label{p11}
	For each $i = 1, \cdots, n -1$, and for any $z \in (a_i, a_{i+1}]$ with $C_i$ given by \eqref{Ci}, we have
	\beqnn
	\mbf{P}_{a_i}\{X_{e_q}\in \d z\} = \frac{qe^{-2\mu_i(a_{i+1} - z)/\sigma_i^2}}{C_i l_i\sigma_i^2 \mbf{E}_{a_{i + 1}}[e^{-q\tau_z}]} \d z
	\eeqnn
	and
	\beqnn
	\mbf{P}_{a_{i+1}}\{X_{e_q}\in \d z\}=  \frac{qe^{-2\mu_i(a_{i+1} - z)/\sigma_i^2}}{C_i l_i\sigma_i^2 \mbf{E}_{a_i}[e^{-q\tau_z}]}  \d z.
	\eeqnn
\end{lem}
\proof
By Corollary \ref{exptau}, for $z \in (a_i, a_{i+1}]$, one obtains that
\beqlb\label{tauz}
	\mbf{E}_{a_i}[e^{-q\tau_z}] = \frac{1}{(1 - c_i^+)e^{\delta^-_i(z - a_i)} + c_i^+ e^{-\delta^+_i(z - a_i)}}
\eeqlb
and
	\beqlb\label{tauz1}
\mbf{E}_{a_{i + 1}}[e^{-q\tau_z}] = \frac{1}{(1 - c_{i+1}^-)e^{\delta^+_i(a_{i+1} - z)} + c_{i+1}^-e^{-\delta^-_i(a_{i+1} - z)}}.
\eeqlb
Notice that
\beqnn
\mbf{P}_{a_i}\{X_{e_q}\in \d z\}\ar=\ar \mbf{P}_{a_i}\{X_{e_q}\in \d z, e_q < \tau_{i+1}\} + \mbf{P}_{a_i}\{X_{e_q}\in \d z, e_q \ge \tau_{i+1}\}\cr
\ar=\ar \mbf{P}_{a_i}\{X_{e_q}\in \d z, e_q < \tau_{i+1}\} + \mbf{P}_{a_i}\{\tau_{i+1} \le e_q\}\mbf{P}_{a_{i+1}}\{X_{e_q} \in \d z\}\cr
\ar=\ar \mbf{P}_{a_i}\{X_{e_q}\in \d z, e_q < \tau_{i+1}\} + \mbf{E}_{a_i}[e^{-q\tau_{i+1}}]\mbf{P}_{a_{i+1}}\{X_{e_q} \in \d z\}
\eeqnn
and
\beqnn
\mbf{P}_{a_{i+1}}\{X_{e_q}\in \d z\}\ar=\ar\mbf{P}_{a_{i+1}}\{X_{e_q}\in \d z, e_q < \tau_{i}\} + \mbf{P}_{a_{i+1}}\{X_{e_q}\in \d z, e_q \ge \tau_{i}\}\cr
\ar=\ar \mbf{P}_{a_{i+1}}\{ X_{e_q}\in \d z, e_q < \tau_{i}\} + \mbf{P}_{a_{i+1}}\{\tau_i \le e_q\}\mbf{P}_{a_i}\{X_{e_q}\in \d z\}\cr
\ar=\ar \mbf{P}_{a_{i+1}}\{X_{e_q}\in \d z, e_q < \tau_{i}\} + \mbf{E}_{a_{i+1}}[e^{-q\tau_i}]\mbf{P}_{a_i}\{X_{e_q}\in \d z\}.
\eeqnn
Then we have
\beqlb\label{eqai}
\ar\ar \mbf{P}_{a_i}\{X_{e_q}\in \d z\} \cr
\ar\ar\quad = \frac{\mbf{P}_{a_i}\{X_{e_q}\in \d z, e_q < \tau_{i+1}\} + \mbf{E}_{a_i}[e^{-q\tau_{i+1}}]  \mbf{P}_{a_{i+1}}\{X_{e_q}\in \d z, e_q < \tau_{i}\}}{1 -  \mbf{E}_{a_i}[e^{-q\tau_{i+1}}] \mbf{E}_{a_{i+1}}[e^{-q\tau_i}]}
\eeqlb
and
\beqlb\label{eqai1}
\ar\ar \mbf{P}_{a_{i+1}}\{X_{e_q}\in \d z\}\cr
\ar\ar\quad = \frac{\mbf{P}_{a_{i+1}}\{X_{e_q}\in \d z, e_q < {\tau}_{i}\} + \mbf{E}_{a_{i+1}}[e^{-q\tau_i}]\mbf{P}_{a_i}\{{X}_{e_q}\in \d z, e_q < {\tau}_{i+1}\} }{1 - \mbf{E}_{a_{i}}[e^{-q\tau_{i+1}}] \mbf{E}_{a_{i+1}}[e^{-q\tau_i}]}.
\eeqlb
Moreover, by \eqref{tauz} and \eqref{tauz1}, one sees that
\beqnn
\ar\ar \frac{1}{\mbf{E}_{a_i}[e^{-q\tau_{i+1}}] \mbf{E}_{a_{i+1}}[e^{-q\tau_{i}}]} - 1 \cr
\ar\ar\quad  = [(1-c^+_{i})e^{\delta_{i}^-(a_{i+1} - a_{i})} + c^+_{i}e^{-\delta_{i}^+(a_{i+1} - a_{i})}] [c^-_{i+1}e^{-\delta_{i}^-(a_{i+1} - a_{i})} + (1- c^-_{i+1})e^{\delta_{i}^+(a_{i + 1} - a_{i})}] - 1\cr
\ar\ar\quad  = (e^{2l_{i}(a_{i+1} - a_{i})} - 1)\left[(1-c^+_{i})(1- c^-_{i+1}) - c^+_{i}c^-_{i+1} e^{-2l_{i}(a_{i+1} - a_{i})}\right]\cr
\ar\ar\quad = 2\sinh((a_{i + 1} - a_i)l_i)\left[(1-c^+_{i})(1- c^-_{i+1})e^{l_i(a_{i + 1} - a_i)} - c^+_{i}c^-_{i+1} e^{-l_{i}(a_{i+1} - a_{i})}\right]
\eeqnn
and
\beqnn
\ar\ar 2\left[\sinh((a_{i+1}-z)l_i) \frac{e^{-\mu_i(a_{i+1} - a_i)/\sigma_i^2}}{ \mbf{E}_{a_{i+1}}[e^{-q\tau_{i}}]} +   \sinh((z-a_{i})l_i)\right] \cr
\ar\ar\quad =2\sinh((a_{i+1} - z)l_i) \left(c_{i+1}^-e^{l_{i}( a_i- a_{i+1})} + (1 - c_{i+1}^-)e^{-l_{i}( a_i - a_{i+1})}\right) + 2\sinh((z - a_i)l_i) \cr
\ar\ar\quad = \left[e^{l_i(a_{i+1} - a_i)} - e^{-l_i(a_{i+1} - a_i)}\right]\left[(1 - c_{i+1}^-)e^{l_i(a_{i+1} - z)} + c_{i+1}^-e^{-l_i(a_{i+1} - z)}\right]\cr
\ar\ar\quad = 2\sinh((a_{i+1} - a_i)l_i) \left[(1 - c_{i+1}^-)e^{l_i(a_{i+1} - z)} + c_{i+1}^-e^{-l_i(a_{i+1} - z)}\right].
\eeqnn
By the above, \eqref{tauz1}, \eqref{eqai}, Lemmas \ref{l0801a} and \ref{l0801b}, we have
\beqnn
\mbf{P}_{a_i}\{X_{e_q}\in \d z\} \ar=\ar \frac{2qe^{\mu_i(z - a_{i})/\sigma_i^2}}{l_i\sigma_i^2} \frac{ \sinh((a_{i+1}-z)l_i) \frac{e^{-\mu_i(a_{i+1} - a_i)/\sigma_i^2}}{ \mbf{E}_{a_{i+1}}[e^{-q\tau_{i}}]} +   \sinh((z-a_{i})l_i) }{\frac{1 }{\mbf{E}_{a_i}[e^{-q\tau_{i+1}}] \mbf{E}_{a_{i+1}}[e^{-q\tau_i}]} - 1}\d z\cr
\ar=\ar  \frac{qe^{\mu_i(z - a_{i})/\sigma_i^2}}{l_i\sigma_i^2} \frac{(1 - c_{i+1}^-)e^{l_i(a_{i+1} - z)} + c_{i+1}^-e^{-l_i(a_{i+1} - z)}}{(1-c^+_{i})(1- c^-_{i+1})e^{l_i(a_{i + 1} - a_i)}  - c^+_{i}c^-_{i+1} e^{- l_{i}(a_{i+1} - a_{i})} }\d z\cr
\ar=\ar  \frac{qe^{-2\mu_i(a_{i+1} - z)/\sigma_i^2}}{C_i l_i\sigma_i^2} \left[(1 - c_{i+1}^-)e^{\delta^+_i(a_{i+1} - z)} + c_{i+1}^-e^{-\delta^-_i(a_{i+1} - z)}\right]\d z\cr
\ar=\ar \frac{qe^{-2\mu_i(a_{i+1} - z)/\sigma_i^2}}{C_i l_i\sigma_i^2\mbf{E}_{a_{i + 1}}[e^{-q\tau_z}]} \d z.
\eeqnn

Similarly, 
\beqnn
\ar\ar  \sinh((z-a_{i})l_i) \frac{e^{\mu_i(a_{i+1} - a_{i})/\sigma_i^2} }{\mbf{E}_{a_{i}}[e^{-q\tau_{i+1}}]} +  \sinh((a_{i+1}-z)l_i) \cr
\ar\ar\quad = \sinh(l_i(a_{i + 1} - a_i))((1 - c_i^+)e^{l_i(z - a_i)} + c_i^+ e^{-l_i(z - a_i)}).
\eeqnn
Combining the above, \eqref{eqai1}, \eqref{tauz} and Lemmas \ref{l0801a} and \ref{l0801b}, we have 
\beqnn
\mbf{P}_{a_{i+1}}\{X_{e_q}\in \d z\} \ar=\ar\frac{2qe^{-\mu_i(a_{i+1} - z)/\sigma_i^2}}{l_i\sigma_i^2 } \frac{ \sinh((z-a_{i})l_i) \frac{e^{\mu_i(a_{i+1} - a_{i})/\sigma_i^2} }{\mbf{E}_{a_{i}}[e^{-q\tau_{i+1}}]} +  \sinh((a_{i+1}-z)l_i)  }{\frac{1}{\mbf{E}_{a_{i}}[e^{-q\tau_{i+1}}] \mbf{E}_{a_{i+1}}[e^{-q\tau_i}]} - 1}\d z\cr
\ar=\ar \frac{qe^{-\mu_i(a_{i+1} - z)/\sigma_i^2}}{l_i\sigma_i^2 } \frac{(1 - c_i^+)e^{l_i(z - a_i)} + c_i^+ e^{-l_i(z - a_i)}} { (1-c^+_{i})(1- c^-_{i+1}) e^{l_i(a_{i + 1} - a_i)} - c^+_{i}c^-_{i+1} e^{- l_{i}(a_{i+1} - a_{i})} }\d z\cr
\ar=\ar \frac{qe^{-2\mu_i(a_{i+1} - z)/\sigma_i^2}}{C_i l_i\sigma_i^2 } \left[(1 - c_i^+)e^{\delta^-_i(z - a_i)} + c_i^+ e^{-\delta^+_i(z - a_i)}\right] \d z\cr
\ar=\ar  \frac{qe^{-2\mu_i(a_{i+1} - z)/\sigma_i^2}}{C_i l_i\sigma_i^2 \mbf{E}_{a_i}[e^{-q\tau_z}]} \d z.
\eeqnn
The proof is complete.
\qed

\begin{prop}\label{p85a}
	For each $i = 1,\cdots, n-1$, and for any $x \le a_i$, $z \in (a_i, a_{i + 1}]$ with $C_i$ given by \eqref{Ci}, we have
		\beqnn
	\mbf{P}_x\{X_{e_q}\in \d z\} = \frac{qe^{-2\mu_i(a_{i+1} - z)/\sigma_i^2}}{C_i l_i\sigma_i^2} \frac{\mbf{E}_x[e^{-q\tau_i}]}{\mbf{E}_{a_{i + 1}}[e^{-q\tau_z}]} \d z.
	\eeqnn
\end{prop}
\proof 
For any $z \in (a_i, a_{i + 1}]$ and $x \le a_i$, we have
\beqnn
\mbf{P}_x\{X_{e_q}\in \d z\} \ar=\ar \mbf{P}_x\{X_{e_q}\in\d z, \tau_i < e_q\} =  \mbf{E}_x[e^{-q\tau_i}]\mbf{P}_{a_i}\{X_{e_q} \in \d z\}.
\eeqnn
The desired result follows by Lemma \ref{p11}.
\qed
\begin{prop}\label{p85b}
		For each $i = 1,\cdots, n-1$, and for any $x \ge a_{i + 1}$, $z \in (a_i, a_{i + 1}]$ with $C_i$ given by \eqref{Ci}, we have
	\beqnn
	\mbf{P}_x\{X_{e_q}\in \d z\} = \frac{qe^{-2\mu_i(a_{i+1} - z)/\sigma_i^2}}{C_i l_i\sigma_i^2 } \frac{\mbf{E}_x[e^{-q\tau_{i + 1}}]}{\mbf{E}_{a_i}[e^{-q\tau_z}]} \d z.
	\eeqnn
\end{prop}
\proof 
For any $z \in (a_i, a_{i + 1}]$ and $x \ge a_{i + 1}$, we have
\beqnn
\mbf{P}_x\{X_{e_q}\in \d z\} \ar=\ar \mbf{P}_x\{X_{e_q}\in\d z, \tau_{i + 1} < e_q\} =  \mbf{E}_x[e^{-q\tau_{i + 1}}]\mbf{P}_{a_{i + 1}}\{X_{e_q} \in \d z\}.
\eeqnn
The result follows by Lemma \ref{p11}.
\qed

\begin{prop}\label{p85c}
	For each $i = 1,\cdots, n-1$, and for any $x, z \in (a_{i}, a_{i+1}]$ with $C_i$ given by \eqref{Ci}, we have
	\beqnn
\mbf{P}_x\{X_{e_q}\in \d z\}  = 	\frac{qe^{-2\mu_i(a_{i+1} - z)/\sigma_i^2}}{C_i l_i\sigma_i^2}\begin{cases}
	\frac{1}{\mbf{E}_{a_i}[e^{-q\tau_x}]\mbf{E}_{a_{i+1}}[e^{-q\tau_z}]} \d z,  & x\le z,\\ 
		 \frac{1}{ \mbf{E}_{a_i}[e^{-q\tau_z}]\mbf{E}_{a_{i+1}}[e^{-q\tau_x}] } \d z,   & z< x.
	\end{cases} 
	\eeqnn
\end{prop}
\proof
	Notice that
		\beqlb\label{eq0809}
\mbf{P}_x\{X_{e_q} \in \d z\}  
	\ar=\ar \mbf{P}_x\{X_{e_q} \in \d z, \tau_i \wedge \tau_{i+1} < e_q, X_{\tau_i \wedge \tau_{i+1}} = a_i\}\cr
	\ar\ar  +   \mbf{P}_x\{X_{e_q} \in \d z, \tau_i \wedge \tau_{i+1} < e_q, X_{\tau_i \wedge \tau_{i+1}} = a_{i + 1}\}\cr
	\ar\ar +  \mbf{P}_x\{X_{e_q} \in \d z, \tau_i \wedge \tau_{i+1} \ge e_q\}\cr
	\ar=\ar \mbf{P}_x\{\tau_i \wedge \tau_{i+1} < e_q, X_{\tau_i \wedge \tau_{i+1}} = a_i\}\mbf{P}_{a_i}\{X_{e_q} \in \d z\}\cr
	\ar\ar +  \mbf{P}_x\{\tau_i \wedge \tau_{i+1} < e_q, X_{\tau_i \wedge \tau_{i+1}} = a_{i + 1}\}\mbf{P}_{a_{i + 1}}\{X_{e_q} \in \d z\}\cr
	\ar\ar +  \mbf{P}_x\{X_{e_q} \in \d z, \tau_i \wedge \tau_{i+1} \ge e_q\}\cr
	\ar=\ar \mbf{P}_x\{\tau_i \wedge \tau_{i+1} < e_q, X^i_{\tau_i \wedge \tau_{i+1}} = a_i\}\mbf{P}_{a_i}\{X_{e_q} \in \d z\}\cr
	\ar\ar +  \mbf{P}_x\{\tau_i \wedge \tau_{i+1} < e_q, X^i_{\tau_i \wedge \tau_{i+1}} = a_{i + 1}\}\mbf{P}_{a_{i + 1}}\{X_{e_q} \in \d z\}\cr
	\ar\ar +  \mbf{P}_x\{X^i_{e_q} \in \d z, \tau_i \wedge \tau_{i+1} \ge e_q\}.
	\eeqlb
The above probabilities are obtained combining \eqref{star4}, \eqref{star5}, Lemma \ref{L01} and Proposition \ref{p11}. In the following we carry out elementary yet tedious calculations to present a simplified result. Notice that
\beqnn
\ar\ar \left[(1 - c_i^+)e^{l_i(z - a_i)} + c_i^+ e^{-l_i(z - a_i)}\right]\cr
\ar\ar\qquad+\left[e^{l_i(a_{i+1}-z)}-e^{-l_i(a_{i+1}-z)}\right] \left[(1-c^+_{i})(1- c^-_{i+1}) e^{l_i(a_{i + 1} - a_i)} - c^+_{i}c^-_{i+1} e^{- l_{i}(a_{i+1} - a_{i})} \right]\cr
\ar\ar\quad = \left[(1 - c_i^+)e^{l_i(a_{i+1} - a_i)} + c_i^+e^{-l_i(a_{i+1} - a_i)} 
\right]\left[(1 - c_{i+1}^-)e^{l_i(a_{i+1} - z)} + c_{i+1}^-e^{-l_i(a_{i+1} - z)}\right]
\eeqnn
and
\beqnn
\ar\ar \left[e^{l_i(a_{i+1}-x)}-e^{-l_i(a_{i+1}-x)}\right] +\left[e^{l_i(x-a_{i})}-e^{-l_i(x-a_i)}\right] \left[(1-c^+_{i}) e^{l_i(a_{i + 1} - a_i)} + c^+_{i}e^{-l_i(a_{i + 1} - a_i)}\right]\cr
\ar\ar\quad = 2\sinh((a_{i + 1} - a_i)l_i)\left[(1 - c_i^+)e^{l_i(x - a_i)} + c_i^+ e^{-l_i(x - a_i)}\right]. 
\eeqnn
Then for the case of $a_i< x\le z\le a_{i+1}$, by the above, \eqref{eq0809}, \eqref{tauz} and \eqref{tauz1}, we have 
	\beqnn
	\ar\ar \mbf{P}_x\{X_{e_q} \in \d z\} \cr
	\ar\ar\quad= \frac{qe^{-2\mu_i(a_{i+1} - z)/\sigma_i^2}}{C_i l_i\sigma_i^2\sinh((a_{i+1} - a_i)l_i)}  \d z\cr
	\ar\ar\qquad\times\bigg[e^{\frac{\mu_i(a_i - x)}{\sigma_i^2}}\sinh((a_{i+1} - x)l_i)[(1 - c_{i+1}^-)e^{\delta^+_i(a_{i+1} - z)} + c_{i+1}^-e^{-\delta^-_i(a_{i+1} - z)}]\cr
	\ar\ar\qquad\quad + e^{\frac{\mu_i(a_{i+1} - x)}{\sigma_i^2}}\sinh((x-a_{i} )l_i)[(1 - c_i^+)e^{\delta^-_i(z - a_i)} + c_i^+ e^{-\delta^+_i(z - a_i)}]\cr
	\ar\ar\qquad\quad +e^{\frac{\mu_i(2a_{i+1}-z - x)}{\sigma_i^2}}2\sinh((x-a_i)l_i)\sinh((a_{i+1}-z)l_i) [(1-c^+_{i})(1- c^-_{i+1}) e^{\delta^-_i(a_{i + 1} - a_i)} \cr
	\ar\ar \qquad\quad - c^+_{i}c^-_{i+1} e^{- \delta^+_{i}(a_{i+1} - a_{i})} ]\bigg]\cr
	\ar\ar\quad= \frac{qe^{-2\mu_i(a_{i+1} - z)/\sigma_i^2} e^{\mu_i(a_i - x+a_{i+1}-z)/\sigma_i^2}}{2C_i l_i\sigma_i^2\sinh((a_{i+1} - a_i)l_i)}   \d z\cr
	\ar\ar\qquad \times\bigg[[e^{l_i(a_{i+1}-x)}-e^{-l_i(a_{i+1}-x)}][(1 - c_{i+1}^-)e^{l_i(a_{i+1} - z)} + c_{i+1}^-e^{-l_i(a_{i+1} - z)}]\cr
	\ar\ar\qquad\quad + [e^{l_i(x-a_{i})}-e^{-l_i(x-a_i)}]\Big[[(1 - c_i^+)e^{l_i(z - a_i)} + c_i^+ e^{-l_i(z - a_i)}]\cr
	\ar\ar\qquad\quad+[e^{l_i(a_{i+1}-z)}-e^{-l_i(a_{i+1}-z)}] [(1-c^+_{i})(1- c^-_{i+1}) e^{l_i(a_{i + 1} - a_i)} - c^+_{i}c^-_{i+1} e^{- l_{i}(a_{i+1} - a_{i})} ]\Big]\bigg]\cr
	\ar\ar\quad=  \frac{qe^{-2\mu_i(a_{i+1} - z)/\sigma_i^2} }{C_i l_i\sigma_i^2}[(1 - c_{i+1}^-)e^{\delta_i^+( a_{i+1}-z)}+c^-_{i+1} e^{- \delta_{i}^-(a_{i + 1} -z)}] \cr
	\ar\ar\qquad\times [(1 - c_i^+)e^{\delta_i^-(x - a_i)} + c_i^+ e^{-\delta_i^+(x - a_i)}] \d z\cr
	\ar\ar\quad= \frac{qe^{-2\mu_i(a_{i+1} - z)/\sigma_i^2}}{C_i l_i\sigma_i^2}\frac{1}{\mbf{E}_{a_i}[e^{-q\tau_x}]\mbf{E}_{a_{i+1}}[e^{-q\tau_z}]} \d z.
	\eeqnn
  
For  $a_i < z< x \le a_{i+1}$, using derivations similar to those above one obtains 
	\beqnn
	\ar\ar \mbf{P}_x\{X_{e_q} \in \d z\} = \frac{qe^{-2\mu_i(a_{i+1} - z)/\sigma_i^2}}{C_i l_i\sigma_i^2} \frac{1}{\mbf{E}_{a_i}[e^{-q\tau_z}]\mbf{E}_{a_{i+1}}[e^{-q\tau_x}]} \d z. 
	\eeqnn
\qed

\subsection{ The case $z > a_n$ }
\begin{lem}\label{p0805a} 
	For any $z > a_n$, we have
	\beqnn
	\mbf{P}_{a_n}\{X_{e_q} \in \d z\}  = \frac{q e^{-\delta_n^-(z - a_n)}}{l_n\sigma_n^2(1-c^+_{n})}\d z.
	\eeqnn
\end{lem}
\proof   
For any $r \in (a_n, z),$ by the Markov property, one sees that
\beqnn
\mbf{P}_{a_n}\{X_{e_q} \in \d z\} \ar=\ar\mbf{P}_{a_n}\{X_{e_q} \in \d z, \tau_r \le e_q\}\cr
\ar=\ar \mbf{P}_{a_n}\{\tau_r \le e_q\}\mbf{P}_r\{X_{e_q} \in \d z\}\cr
\ar=\ar \mbf{E}_{a_n}[e^{-q\tau_r}]\left[\mbf{P}_r\{X_{e_q} \in \d z, \tau_n > e_q\} + \mbf{P}_r\{X_{e_q} \in \d z, \tau_n \le e_q\}\right]\cr
\ar=\ar \mbf{E}_{a_n}[e^{-q\tau_r}]\left[\mbf{P}_r\{X^n_{e_q} \in \d z, \tau_n > e_q\} + \mbf{P}_r\{\tau_n \le e_q\}\mbf{P}_{a_n}\{X_{e_q} \in \d z\}\right]\cr
\ar=\ar \mbf{E}_{a_n}[e^{-q\tau_r}]\left[\mbf{P}_r\{X^n_{e_q} \in \d z, \tau_n > e_q\} + \mbf{E}_r[e^{-q\tau_n}]\mbf{P}_{a_n}\{X_{e_q} \in \d z\}\right],
\eeqnn
which implies that
\beqlb\label{2}
\mbf{P}_{a_n}\{X_{e_q} \in \d z\} = \frac{\mbf{E}_{a_n}[e^{-q\tau_r}]\mbf{P}_r\{X^n_{e_q} \in \d z, \tau_n > e_q\}}{1 - \mbf{E}_{a_n}[e^{-q\tau_r}]\mbf{E}_r[e^{-q\tau_n}]}.
\eeqlb
Notice that $\tau_n = \inf\{t \ge 0: X_t^n = a_n\}$ when $X_0 = X_0^n = r$. It follows from \eqref{1} that 
\beqnn
\mbf{P}_r\{X^n_{e_q} \in \d z, \tau_n > e_q\} = \frac{q}{l_n\sigma_n^2}\left(e^{\delta_n^-(r - a_n)} - e^{-\delta_n^+(r - a_n)}\right)e^{-\delta_n^-(z - a_n)}\d z.
\eeqnn
Moreover, by Corollary \ref{exptau}, we have
\beqlb\label{eqanr}
\mbf{E}_{a_n}[e^{-q\tau_r}] = \frac{g_q^+(a_n)}{g_q^+(r)} = \frac{1}{(1-c^+_{n})e^{\delta_{n}^-(r-a_n)}+c^+_{n}e^{-\delta_{n}^+(r-a_n)}}
\eeqlb
and
\beqnn
\mbf{E}_r[e^{-q\tau_n}] = \frac{g_q^-(r)}{g_q^-(a_n)} = e^{-\delta_n^+(r - a_n)}.
\eeqnn
One sees that
\beqlb\label{eq0805a}
\mbf{P}_{a_n}\{X_{e_q} \in \d z\} \ar=\ar \frac{\mbf{P}_r\{X^n_{e_q} \in \d z, \tau_n > e_q\}}{\frac{1}{\mbf{E}_{a_n}[e^{-q\tau_r}]} - \mbf{E}_r[e^{-q\tau_n}]}\cr
\ar=\ar \frac{\frac{q}{l_n\sigma_n^2}(e^{\delta_n^-(r - a_n)} - e^{-\delta_n^+(r - a_n)})e^{-\delta_n^-(z - a_n)}\d z}{(1-c^+_{n})e^{\delta_{n}^-(r-a_n)}+c^+_{n}e^{-\delta_{n}^+(r-a_n)} -  e^{-\delta_n^+(r - a_n)}}\cr
\ar=\ar \frac{q e^{-\delta_n^-(z - a_n)}}{l_n\sigma_n^2(1-c^+_{n})}\d z.
\eeqlb
The proof is completed.
\qed

\begin{prop}\label{p0717a}
	For any $x, z > a_n$, we have
	\beqnn
	\mbf{P}_{x}\{X_{e_q}\in \d z\}
	\ar=\ar \begin{cases}
		\frac{qe^{-\delta_n^-(z - a_n)}/\mbf{E}_{a_n}[e^{-q\tau_x}]}{l_n\sigma_n^2 (1 - c_n^+)} \d z, & z \ge x > a_n,\\
		e^{2\mu_n(z - x)/\sigma_n^2}\frac{qe^{-\delta^-_n(x - a_n)}/\mbf{E}_{a_n}[e^{-q\tau_z}]}{l_n\sigma_n^2(1 - c_n^+)}\d z, & x > z > a_n.
	\end{cases}
	\eeqnn	
\end{prop}

\proof
For any $x, z \ge a_n,$ by the Markov property, we have
\beqlb\label{eq311}
\mbf{P}_{x}\{X_{e_q}\in \d z\}\ar=\ar \mbf{P}_{x}\{X_{e_q}\in \d z, e_q \le \tau_n\} + \mbf{P}_x\{X_{e_q}\in \d z, e_q > \tau_n\}\cr
\ar=\ar \mbf{P}_{x}\{X_{e_q}^n\in \d z, e_q \le \tau_n^n\} + \mbf{P}_x\{\tau_n < e_q\}\mbf{P}_{a_n}\{X_{e_q}\in \d z\}\cr
\ar=\ar  \mbf{P}_{x}\{X_{e_q}^n\in \d z, e_q \le \tau_n^n\} + \mbf{E}_x[e^{-q\tau_n}]\mbf{P}_{a_n}\{X_{e_q}\in \d z\}.
\eeqlb
By \eqref{1}, recall that
\beqlb\label{eq311a}
\ar\ar \mbf{P}_{x}\{X_{e_q}^n\in \d z, e_q \le \tau_n^n\}\cr
\ar\ar\quad = \begin{cases}
	\frac{q}{l_n\sigma_n^2}(e^{\delta_n^-(x - a_n)} - e^{-\delta_n^+(x - a_n)})e^{-\delta_n^-(z - a_n)}\d z,  & z \ge x > a_n,\\
	\frac{q}{l_n\sigma_n^2}(e^{\delta_n^+(z - a_n)} - e^{-\delta_n^-(z - a_n)})e^{-\delta_n^+(x - a_n)}\d z,  & x > z > a_n.
\end{cases}
\eeqlb
Moreover, by Corollary \ref{exptau}, one obtains that
\beqlb\label{eq311b}
\mbf{E}_x[e^{-q\tau_n}] = \frac{g_q^-(x)}{g_q^-(a_n)} = e^{-\delta_n^+(x - a_n)}
\eeqlb
and
\beqlb\label{eq37a}
\mbf{E}_{a_n}[e^{-q\tau_x}] \ar=\ar \frac{1}{(1 - c_n^+)e^{\delta_n^-(x - a_n)} + c_n^+ e^{-\delta_n^+(x - a_n)}}\cr
\ar=\ar \frac{e^{\mu_n(x - a_n)/\sigma_n^2}}{(1 - c_n^+)e^{l_n(x - a_n)} + c_n^+ e^{-l_n(x - a_n)}}.
\eeqlb
It follows from \eqref{eq311}-\eqref{eq311b} and Lemma \ref{p0805a} that 
\beqnn
\mbf{P}_{x}\{X_{e_q}\in \d z\}
\ar=\ar\begin{cases}
	\frac{qe^{-\delta_n^-(z - a_n)}}{l_n\sigma_n^2}\left(e^{\delta_n^-(x - a_n)} - e^{-\delta_n^+(x - a_n)} + \frac{e^{-\delta_n^+(x - a_n)}}{1 - c_n^+}\right)\d z,  & z \ge x > a_n,\\
	\frac{qe^{-\delta_n^+(x - a_n)}}{l_n\sigma_n^2}\left(e^{\delta_n^+(z - a_n)} - e^{-\delta_n^-(z - a_n)} + \frac{e^{-\delta_n^-(z - a_n)}}{1 - c_n^+}\right)\d z,  & x > z > a_n
\end{cases}\cr
\ar=\ar \begin{cases}
	e^{\mu_n(z - x)/\sigma_n^2}\frac{qe^{-l_n(z - a_n)}}{l_n\sigma_n^2}\left(e^{l_n(x - a_n)} + \frac{c_n^+}{1 - c_n^+}e^{-l_n(x - a_n)}\right)\d z,  & z \ge x > a_n,\\
	e^{\mu_n(z - x)/\sigma_n^2}\frac{qe^{-l_n(x - a_n)}}{l_n\sigma_n^2}\left(e^{l_n(z - a_n)}  + \frac{c_n^+}{1 - c_n^+} e^{-l_n(z - a_n)}\right)\d z,  & x > z > a_n.
\end{cases}
\eeqnn
The proof is completed by the above and \eqref{eq37a}.
\qed

\begin{prop}\label{p0805b}
	For any $z > a_n \ge x$, we have
	\beqnn
	\mbf{P}_x\{X_{e_q} \in \d z\} =  \mbf{E}_x[e^{-q\tau_n}]\frac{q e^{-\delta_n^-(z - a_n)}}{l_n\sigma_n^2(1-c^+_{n})}\d z.
	\eeqnn
\end{prop}
\proof 
For any $z > a_n \ge x$, by Corollary \ref{exptau} we have
\beqnn
\mbf{P}_x\{X_{e_q} \in \d z\} \ar=\ar \mbf{P}_x\{X_{e_q} \in \d z, \tau_n < e_q\} = \mbf{P}_x\{\tau_n < e_q\}\mbf{P}_{a_n}\{X_{e_q}\in \d z\}\cr
\ar=\ar \mbf{E}_x[e^{-q\tau_n}] \mbf{P}_{a_n}\{X_{e_q}\in \d z\}.
\eeqnn
The result follows by Lemma \ref{p0805a}.
\qed

{\bf Proof of Theorem \ref{tpm}} 
(i) For $z \le a_1$, 
the desired result follows from Propositions \ref{p0717b} and \ref{p0805d} together with the  Corollary \ref{exptau}.

(ii) For $z \in (a_i, a_{i+1}]$ for $i = 1, \cdots, n-1$,
the desired results follow from sequential applications  of Propositions \ref{p85a}--\ref{p85c} and Corollary \ref{exptau}.

(iii) For $z > a_n$,
the result follows from a combination of Propositions \ref{p0717a} and \ref{p0805b} in conjunction with Corollary \ref{exptau}.
\qed

\section{Proof of Corollary \ref{main1}}\label{SI}
 
 In this section, we present the proof of Corollary \ref{main1} under the assumption of $\mu_0 > 0$ and $\mu_n < 0$.

	\begin{lem}\label{l0729a}
	Assume that $\mu_0 > 0$ and $\mu_i \neq 0$ for any $i = 1, \cdots, n$. Then we have
	\beqnn
	\lim_{q \rightarrow 0}c_i^+ = \begin{cases}
		1, & \mu_i < 0,\\
		0, & \mu_i > 0.
	\end{cases}
	\eeqnn
\end{lem}

\proof
By \eqref{delta}, one sees that 
\beqnn
\lim_{q \rightarrow 0}\delta_i^- = \frac{|\mu_i| - \mu_i}{\sigma_i^2}, \quad \lim_{q \rightarrow 0}\delta_i^+ = \frac{|\mu_i| + \mu_i}{\sigma_i^2}, \quad \lim_{q \rightarrow 0}l_i = \frac{|\mu_i|}{\sigma_i^2}
\eeqnn
for any $i = 1, \cdots, n.$ Moreover,
\beqnn
\lim_{q \rightarrow 0} \frac{\delta_i^- - \delta_{i-1}^-}{\delta_i^+ + \delta_i^-} =\frac12\left[1 - \frac{\mu_i}{|\mu_i|} - \frac{\sigma_i^2}{|\mu_i|}\frac{|\mu_{i-1}| - \mu_{i - 1}}{\sigma_{i-1}^2}\right]. 
\eeqnn
Since $\mu_0 > 0$,  by \eqref{c^+_i} we have
\beqnn
\lim_{q \rightarrow 0}c_1^+ \ar=\ar \lim_{q \rightarrow 0}\frac{\delta_1^- - \delta_0^-}{\delta_1^+ + \delta_1^-} = \frac12\left[1 - \frac{\mu_1}{|\mu_1|} - \frac{\sigma_1^2}{|\mu_1|}\frac{|\mu_{0}| - \mu_{0}}{\sigma_{0}^2}\right] \cr
\ar=\ar \frac12\left[1 - \frac{\mu_1}{|\mu_1|}\right] = \begin{cases}
	0, & \mu_1 > 0,\\
	1, & \mu_1 < 0.
\end{cases}
\eeqnn
Similarly,  assume that $\lim_{q \rightarrow 0}c_{i - 1}^+ = 1$ if $\mu_{i - 1} < 0$ and $\lim_{q \rightarrow 0}c_{i - 1}^+ = 0$ if $\mu_{i - 1} > 0.$ One obtains that
\beqnn
\lim_{q \rightarrow 0}\frac{c^+_{i - 1}}{(1-c^+_{i - 1})e^{2l_{i - 1}(a_{i} - a_{i - 1})} + c^+_{i - 1}} = \begin{cases}
	0, & \mu_{i-1}> 0,\\
	1,& \mu_{i-1} < 0,
\end{cases}
\eeqnn
which implies 
\beqnn
\lim_{q \rightarrow 0} c_i^+ \ar=\ar \lim_{q \rightarrow 0}\frac{1}{2l_i}\left[(\delta_i^- - \delta_{i-1}^-) +  \frac{2l_{i-1}c^+_{i-1}}{(1-c^+_{i-1})e^{2l_{i - 1}(a_{i} - a_{i-1})} + c^+_{i-1}}\right]\cr
\ar=\ar \frac{\sigma_i^2}{2|\mu_i|}\left[\frac{|\mu_i| - \mu_i}{\sigma_i^2} - \frac{|\mu_{i - 1}| - \mu_{i - 1}}{\sigma_{i - 1}^2} + \frac{2|\mu_{i - 1}|}{\sigma_{i - 1}^2}\lim_{q \rightarrow 0}\frac{c^+_{i - 1}}{(1-c^+_{i - 1})e^{2l_{i - 1}(a_{i} - a_{i - 1})} + c^+_{i - 1}}\right]\cr
\ar=\ar \frac{\sigma_i^2}{2|\mu_i|} \frac{|\mu_i| - \mu_i}{\sigma_i^2} =  \begin{cases}
	0, & \mu_i > 0,\\
	1, & \mu_i < 0.
\end{cases}
\eeqnn
Then the result follows by induction. 
\qed

If $\mu_i < 0$, then
\beqlb\label{eq0603a}
\lim_{q \rightarrow 0}\frac{\delta_i^+}{q} = \lim_{q \rightarrow 0}\frac{\sqrt{2q\sigma_i^2 + \mu_i^2} + \mu_i}{q\sigma_i^2} = \lim_{q \rightarrow 0}\frac{2}{\sqrt{2q\sigma_i^2 + \mu_i^2} - \mu_i} = -\frac{1}{\mu_i}.
\eeqlb
On the other hand, if $\mu_i > 0$, we have 
\beqlb\label{eq0603b}
\lim_{q \rightarrow 0}\frac{\delta_i^-}{q} = \lim_{q \rightarrow 0}\frac{\sqrt{2q\sigma_i^2 + \mu_i^2} - \mu_i}{q\sigma_i^2} = \lim_{q \rightarrow 0}\frac{2}{\sqrt{2q\sigma_i^2 + \mu_i^2} + \mu_i} = \frac{1}{\mu_i}.
\eeqlb
 Before stating the result, we define a sequence $\{F_i: i = 1, 2, \cdots, n\}$  satisfying 
\beqlb\label{F}
\begin{cases}
	F_i = \frac{1}{\mu_{i - 1}}{\bf 1}_{\{\mu_{i-1} \neq 0\}} + \frac{2(a_i - a_{i - 1})}{\sigma_{i - 1}^2}{\bf 1}_{\{\mu_{i - 1} = 0\}} - \frac{1}{\mu_i}{\bf 1}_{\{\mu_i \neq 0\}}  + e^{-\frac{2\mu_{i - 1}(a_i - a_{i - 1})}{\sigma_{i - 1}^2}}F_{i - 1},\\
	F_1 = \frac{1}{\mu_0} - \frac{1}{\mu_1}{\bf 1}_{\{\mu_1 \neq 0\}}.
\end{cases}
\eeqlb
Then we can get the following result.
\begin{lem}\label{l0730a}
	If $\mu_0 > 0$ and $\mu_i \neq 0$ for  $i = 1, \cdots, n$, for $F_i$ given by \eqref{F} we have
	\beqnn
	\lim_{q \rightarrow 0}\frac{2l_i[(1 - c_i^+){\bf 1}_{\{\mu_i < 0\}} - c_i^+{\bf 1}_{\{\mu_i > 0\}}]}{q} = F_i.
	\eeqnn
\end{lem}

\proof
If $\mu_i < 0$ and $\mu_{i - 1} < 0$, by Lemma \ref{l0729a}, we have $\lim_{q \rightarrow 0}c_i^+ = \lim_{q \rightarrow 0}c_{i - 1}^+ = 1.$ Then by \eqref{eq0603a}, one sees that
\beqnn
\lim_{q \rightarrow 0}\frac{2l_i(1 - c_i^+)}{q} 
\ar=\ar  \lim_{q \rightarrow 0} \left[\frac{\delta_i^+}{q} + \frac{1}{q}\left(\delta_{i - 1}^- - \frac{2l_{i - 1}c^+_{i-1}}{(1-c^+_{i-1})e^{2l_{i - 1}(a_{i} - a_{i-1})} + c^+_{i-1}} \right)\right] \cr
\ar=\ar  \lim_{q \rightarrow 0} \left[\frac{\delta_i^+}{q} + \frac{1}{q} \frac{(1-c^+_{i-1})\delta_{i - 1}^- e^{2l_{i - 1}(a_{i} - a_{i-1})} - \delta_{i - 1}^+c^+_{i-1}}{(1-c^+_{i-1})e^{2l_{i - 1}(a_{i} - a_{i-1})} + c^+_{i-1}} \right] \cr
\ar=\ar \frac{1}{\mu_{i - 1}}-\frac{1}{\mu_i} + e^{-2\mu_{i - 1}(a_i - a_{i - 1})/\sigma_{i - 1}^2}\lim_{q \rightarrow 0}\frac{2l_{i - 1}(1 - c_{i - 1}^+)}{q}.
\eeqnn
If $\mu_i < 0$ and $\mu_{i - 1} > 0$, then we have $\lim_{q \rightarrow 0}c_i^+ = 1$ and $\lim_{q \rightarrow 0}c_{i - 1}^+ = 0$ by Lemma \ref{l0729a}. It follows from \eqref{eq0603a} and \eqref{eq0603b} that
\beqnn
\lim_{q \rightarrow 0}\frac{2l_i(1 - c_i^+)}{q} 
\ar=\ar \lim_{q \rightarrow 0} \left[\frac{\delta_i^+ + \delta_{i-1}^-}{q} - \frac{1}{q}\frac{2l_{i - 1}c^+_{i-1}}{(1-c^+_{i-1})e^{2l_{i - 1}(a_{i} - a_{i-1})} + c^+_{i-1}} \right]\cr
\ar=\ar \frac{1}{\mu_{i - 1}} - \frac{1}{\mu_i} - e^{-2\mu_{i - 1}(a_i - a_{i - 1})/\sigma_{i - 1}^2}\lim_{q \rightarrow 0}\frac{2l_{i - 1}c_{i - 1}^+}{q}.
\eeqnn
If $\mu_i > 0$ and $\mu_{i - 1}<0$, then we have $\lim_{q \rightarrow 0}c_i^+ = 0$ and $\lim_{q \rightarrow 0}c_{i - 1}^+ = 1$ by Lemma \ref{l0729a}. One obtains from \eqref{eq0603a} and \eqref{eq0603b} that
\beqnn
\lim_{q \rightarrow 0}\frac{-2l_i c_i^+}{q} \ar=\ar-\lim_{q \rightarrow 0}\frac{1}{q}\left[ \delta_i^- - \delta_{i-1}^- + \frac{2l_{i - 1}c^+_{i-1}}{(1-c^+_{i-1})e^{2l_{i - 1} (a_{i} - a_{i-1})} + c^+_{i-1}}\right]\cr
\ar=\ar -\lim_{q \rightarrow 0}\left[\frac{\delta_i^-}{q} - \frac{1}{q}\frac{(1-c^+_{i-1})\delta_{i - 1}^-e^{2l_{i - 1}(a_{i} - a_{i-1})} - \delta_{i - 1}^+c_{i - 1}^+}{(1-c^+_{i-1})e^{2l_{i - 1} (a_{i} - a_{i-1})} + c^+_{i-1}}\right]\cr
\ar=\ar \frac{1}{\mu_{i - 1}} - \frac{1}{\mu_i} + e^{-2\mu_{i - 1}(a_i - a_{i - 1})/\sigma_{i - 1}^2}\lim_{q \rightarrow 0}\frac{2l_{i - 1}(1 - c_{i - 1}^+)}{q}.
\eeqnn
If $\mu_i > 0$ and $\mu_{i - 1} > 0$, then by Lemma \ref{l0729a}, one sees that $\lim_{q \rightarrow 0}c_i^+ = \lim_{q \rightarrow 0}c_{i - 1}^+ = 0$. It follows from  \eqref{eq0603b} that
\beqnn
\lim_{q \rightarrow 0}\frac{-2l_i c_i^+}{q} \ar=\ar - \lim_{q \rightarrow 0}\frac{1}{q}\left[ \delta_i^- - \delta_{i-1}^- + \frac{2l_{i - 1}c^+_{i-1}}{(1-c^+_{i-1})e^{2l_{i - 1} (a_{i} - a_{i-1})} + c^+_{i-1}}\right]\cr
\ar=\ar \frac{1}{\mu_{i - 1}} - \frac{1}{\mu_i} + e^{-2\mu_{i - 1}(a_i - a_{i - 1})/\sigma_{i - 1}^2}\lim_{q \rightarrow 0}\frac{-2l_{i - 1}c_{i - 1}^+}{q}.
\eeqnn
Moreover, we have $\lim_{q \rightarrow 0}c_1^+ = 0$ if $\mu_1 > 0$. Then
\beqnn
\lim_{q \rightarrow 0}\frac{-2l_1c_1^+}{q} \ar=\ar -\lim_{q \rightarrow 0}\frac{\delta_1^- - \delta_0^-}{q} = -\lim_{q \rightarrow 0}\frac{\frac{\sqrt{2q\sigma_1^2 + \mu_1^2} - \mu_1}{\sigma_1^2} - \frac{\sqrt{2q\sigma_0^2 + \mu_0^2} - \mu_0}{\sigma_0^2}}{q} \cr
\ar=\ar -\lim_{q \rightarrow 0}\frac{\frac{2q}{\sqrt{2q\sigma_1^2 + \mu_1^2} + \mu_1} - \frac{2q}{\sqrt{2q\sigma_0^2 + \mu_0^2} + \mu_0}}{q} = \frac{1}{\mu_0} - \frac{1}{\mu_1}.
\eeqnn
And $\lim_{q \rightarrow 0}c_1^+ = 1$ if $\mu_1 < 0$. Then
\beqnn
\lim_{q \rightarrow 0}\frac{2l_1(1 - c_1^+)}{q} \ar=\ar \lim_{q \rightarrow 0}\frac{\delta_1^+ + \delta_0^-}{q} = \lim_{q \rightarrow 0}\frac{\frac{\sqrt{2q\sigma_1^2 + \mu_1^2} + \mu_1}{\sigma_1^2} + \frac{\sqrt{2q\sigma_0^2 + \mu_0^2} - \mu_0}{\sigma_0^2}}{q} \cr
\ar=\ar \lim_{q \rightarrow 0}\frac{\frac{2q}{ \sqrt{2q\sigma_1^2 + \mu_1^2} - \mu_1} + \frac{2q}{\sqrt{2q\sigma_0^2 + \mu_0^2} + \mu_0}}{q} = \frac{1}{\mu_0} - \frac{1}{\mu_1}.
\eeqnn
The desired result follows.
\qed

	\begin{cor}\label{c0724}
		If $\mu_0 > 0$, then for any $i = 1, \cdots, n$, we have
		\beqnn
		\lim_{q \rightarrow 0}c_i^+ = \begin{cases}
			1, & \mu_i < 0,\\
			1/2, & \mu_i = 0,\\
			0, & \mu_i > 0.
		\end{cases}
		\eeqnn
	\end{cor}
\proof 
For the case of $\mu_1 = 0$, we have
\beqnn
\lim_{q \rightarrow 0}c_1^+ =  \lim_{q \rightarrow 0}\frac{\delta_1^- - \delta_0^-}{\delta_1^+ + \delta_1^-} = \frac{1}{2}.
\eeqnn 
 Notice that $\delta_i^- = \delta_i^+ = l_i = \frac{\sqrt{2q} }{|\sigma_i|}$ when $\mu_i = 0.$ For the case of $\mu_i = 0$ and $\mu_{i-1} \neq 0$,   by \eqref{eq0603a}, \eqref{eq0603b} and Lemma \ref{l0730a}, one can check that
	\beqnn
	\lim_{q \rightarrow 0} c_i^+ \ar=\ar \lim_{q \rightarrow 0}  \frac{\delta_i^- - \delta_{i-1}^-}{\delta_i^+ + \delta_i^-} + \frac{\delta_{i-1}^+ + \delta_{i-1}^-}{\delta_i^+ + \delta_i^-} \frac{c^+_{i-1}}{(1-c^+_{i-1})e^{(\delta_{i-1}^- + \delta_{i - 1}^+)(a_{i} - a_{i-1})} + c^+_{i-1}}\cr
	\ar=\ar \frac{1}{2} + \frac{|\sigma_i|}{2}\lim_{q \rightarrow 0}\frac{1}{\sqrt{2q}}\frac{\delta_{i-1}^+ c^+_{i-1} - \delta_{i-1}^- (1-c^+_{i-1})e^{2l_{i-1} (a_{i} - a_{i-1})}  }{(1-c^+_{i-1})e^{2l_{i-1} (a_{i} - a_{i-1})} + c^+_{i-1}} \cr
	\ar=\ar  \frac{1}{2} + \frac{|\sigma_i|}{2}\lim_{q \rightarrow 0}\frac{q}{\sqrt{2q}}\frac{\frac{\delta_{i-1}^+c^+_{i-1}}{q} - \frac{\delta_{i-1}^- (1-c^+_{i-1})}{q}e^{2l_{i-1} (a_{i} - a_{i-1})}  }{(1-c^+_{i-1})e^{2l_{i-1} (a_{i} - a_{i-1})} + c^+_{i-1}} \cr
	\ar=\ar \frac{1}{2}.
	\eeqnn 
	 
For the case of $\mu_i = 0$ and $\mu_{i-1} = 0$, we have $\delta_{i-1}^- = \delta_{i-1}^+ = l_{i-1} = \frac{\sqrt{2q}}{|\sigma_{i-1}|}.$ Then
	\beqnn
	\lim_{q \rightarrow 0} c_i^+ = \frac{1}{2} - \frac{|\sigma_i|}{2|\sigma_{i-1}|} + \frac{|\sigma_i|}{2|\sigma_{i-1}|} = \frac{1}{2}.
	\eeqnn 
	
On the other hand, for the case of $\mu_{i-1} = 0$ and $\mu_i \neq 0$,	similar to the proof of Lemma \ref{l0729a}, one can check that 
		\beqnn
	\lim_{q \rightarrow 0}c_i^+ = \begin{cases}
		1, & \mu_i < 0,\\
		0, & \mu_i > 0.
	\end{cases}
	\eeqnn
 Then the result follows from induction method and Lemma \ref{l0729a}.
	\qed
	
	\begin{cor}
			If $\mu_0 > 0$, then for $i = 1, \cdots, n$ we have
		\beqlb\label{eq0725}
		\lim_{q \rightarrow 0}\frac{2l_i[(1 - c_i^+){\bf 1}_{\{\mu_i < 0\}} + (\frac12 - c_i^+){\bf 1}_{\{\mu_i = 0\}} - c_i^+{\bf 1}_{\{\mu_i > 0\}}]}{q} = F_i,
		\eeqlb
		where $F_i$ is given by \eqref{F}.
	\end{cor}
\proof 
For the case of  $\mu_1 = 0$, we have $\delta_1^+ = \delta_1^- = l_1 = \frac{\sqrt{2q}}{|\sigma_1|}$. Then by \eqref{eq0603b}, one obtains that
\beqnn
\lim_{q \rightarrow 0}\frac{2l_1(\frac12 - c_1^+)}{q} \ar=\ar \lim_{q \rightarrow 0} \frac{\delta_0^-}{q} = \frac{1}{\mu_0}.
\eeqnn 

Now we assume \eqref{eq0725} holds for $i-1$. For the case of $\mu_i = 0$ and $\mu_{i-1} \neq 0$, one sees that
	\beqnn
	\lim_{q \rightarrow 0} \frac{2l_i(\frac{1}{2} - c_i^+)}{q} \ar=\ar  \lim_{q \rightarrow 0} \frac{1}{q} \frac{ \delta_{i-1}^- (1-c^+_{i-1})e^{2l_{i-1} (a_{i} - a_{i-1})} - \delta_{i-1}^+ c^+_{i-1}}{(1-c^+_{i-1})e^{2l_{i-1} (a_{i} - a_{i-1})} + c^+_{i-1}} \cr
	\ar=\ar \frac{1}{\mu_{i-1}} + e^{-2\mu_{i-1}(a_i - a_{i-1})/\sigma_{i-1}^2}F_{i-1}.
	\eeqnn 
For the case of $\mu_i = \mu_{i - 1} = 0$, we have
	\beqnn
	\lim_{q \rightarrow 0} \frac{2l_i(\frac{1}{2} - c_i^+)}{q} \ar=\ar  \lim_{q \rightarrow 0} \frac{1}{q} \frac{\frac{\sqrt{2q}}{|\sigma_{i-1}|} (1-c^+_{i-1})e^{2\frac{\sqrt{2q}}{|\sigma_{i-1}|} (a_{i} - a_{i-1})} - \frac{\sqrt{2q}}{|\sigma_{i-1}|} c^+_{i-1}}{(1-c^+_{i-1})e^{2\frac{\sqrt{2q}}{|\sigma_{i-1}|} (a_{i} - a_{i-1})} + c^+_{i-1}} \cr
	\ar=\ar \lim_{q \rightarrow 0}\frac{2l_{i-1}(\frac{1}{2} - c_{i-1}^+)}{q} + \lim_{q \rightarrow 0}\frac{\frac{\sqrt{2q}}{|\sigma_{i-1}|}c_{i-1}^+ (e^{2\frac{\sqrt{2q}}{|\sigma_{i-1}|}(a_i - a_{i-1})} - 1)}{q} \cr
	\ar=\ar  \lim_{q \rightarrow 0}\frac{2l_{i-1}(\frac{1}{2} - c_{i-1}^+)}{q} + \frac{2(a_i - a_{i-1})}{\sigma_{i-1}^2}\cr
	\ar=\ar F_{i-1} + \frac{2(a_i - a_{i-1})}{\sigma_{i-1}^2}.
	\eeqnn 
	On the other hand, for the case of $\mu_{i-1} = 0$ and $\mu_i \neq 0$, similar to the proof of Lemma \ref{l0730a} and the above, we have
		\beqnn
	\lim_{q \rightarrow 0} \frac{2l_i[(1 - c_i^+){\bf 1}_{\{\mu_i < 0\}} - c_i^+{\bf 1}_{\{\mu_i > 0\}}}{q} \ar=\ar  \frac{2(a_i - a_{i-1})}{\sigma_{i-1}^2} - \frac{1}{\mu_i} + F_{i-1}.
	\eeqnn
	Then the result follows from induction method and Lemma \ref{l0730a}.
	\qed

Let 
\beqlb\label{F1}
\begin{cases}
	\overline{F}_i =  \frac{1}{\mu_{i - 1}}{\bf 1}_{\{\mu_{i-1} \neq 0\}} - \frac{1}{\mu_i}{\bf 1}_{\{\mu_i \neq 0\}} + \frac{2(a_{i+1} - a_i)}{\sigma_i^2}{\bf 1}_{\{\mu_i = 0\}} + e^{2\mu_{i }(a_{i+1} - a_{i })/\sigma_{i }^2}\overline{F}_{i + 1},\\
	\overline{F}_n =  \frac{1}{\mu_{n-1}}{\bf 1}_{\{\mu_{n - 1} \neq 0\}} - \frac{1}{\mu_n}.
\end{cases}
\eeqlb
It is easy to check $\overline{F}_1$ given by \eqref{overlineF1} satisfies the above \eqref{F1}.
By the demonstration methodology from Lemma \ref{l0730a} and Corollary \ref{c0724}, we establish the subsequent results, presented here without proof.

	\begin{lem}\label{l0729b}
	If $\mu_n < 0$, then for any $i = 1, \cdots, n$, we have
	\beqnn
	\lim_{q \rightarrow 0}c_i^- = \begin{cases}
		0, & \mu_{i- 1} < 0,\\
		1/2, & \mu_{i - 1} = 0,\\
		1, & \mu_{i- 1} > 0.
	\end{cases}
	\eeqnn
\end{lem}

\begin{lem}\label{l0730b}
	If $\mu_n < 0$, then for $i = 1, \cdots, n$ and $\overline{F}_i$  given by \eqref{F1}, we have
	\beqnn
	\lim_{q \rightarrow 0}\frac{2l_{i-1}\left[(1 - c_i^-){\bf 1}_{\{\mu_{i - 1} > 0\}} + (\frac12 - c_i^-){\bf 1}_{\{\mu_{i - 1} = 0\}} - c_i^-{\bf 1}_{\{\mu_{i - 1} < 0\}}\right]}{q} = \overline{F}_i.
	\eeqnn

\end{lem}

\begin{prop}\label{p312a}
		If $\mu_0 > 0$ and $\mu_n < 0$, then for $z \le a_1, x \in \mathbb{R}$ and $\overline{F}_1$ given by \eqref{overlineF1}, we have
		\beqnn
		\lim_{q \rightarrow 0}\mbf{P}_x\{X_{e_q} \in \d z\} = \frac{2e^{-2\mu_0(a_1 - z)/ \sigma_0^2}}{\sigma_0^2 \overline{F}_{1}}\d z.
		\eeqnn

\end{prop}

\proof 
 Under the assumption of $\mu_0 > 0$ and $\mu_n < 0$, one sees that $\lim_{q \rightarrow 0}\delta_0^+ = 2\mu_0/\sigma_0^2$, $\lim_{q \rightarrow 0}\delta_0^- = 0$ and $\lim_{q \rightarrow 0}c_1^- = 1$ by \eqref{delta} and Lemma \ref{l0729b}.
 
 For the case of $a_1 \ge z \ge x,$ by Corollary \ref{exptau}, we have 
 \beqnn
 \lim_{q \rightarrow 0}\frac{q_q^-(z)}{b_1^-} = \lim_{q \rightarrow 0}\frac{g_q^-(z)}{g_q^-(a_1)} = \lim_{q \rightarrow 0}\frac{1}{\mbf{E}_{a_1}[e^{-q\tau_z}]} = 1.
 \eeqnn 
 Then by Lemma \ref{l0730b} and Theorem \ref{tpm} we have
 \beqnn
 	\lim_{q \rightarrow 0}\mbf{P}_x\{X_{e_q} \in \d z\} \ar=\ar \lim_{q \rightarrow 0} \frac{q e^{-\delta_0^+(a_1 - z)}}{l_0\sigma_0^2 b_1^- (1-c^-_1)}	g_q^-(z) e^{-\delta_0^-(z - x)} \d z\cr
 	\ar=\ar  \frac{2e^{-2\mu_0(a_1 - z)/ \sigma_0^2}}{\sigma_0^2 \overline{F}_{1}}\d z.
 \eeqnn 
 
 For $x \ge z$, notice that $\lim_{q \rightarrow 0}\frac{q_q^-(x)}{b_1^-} =\lim_{q \rightarrow 0} \frac{g_q^-(x)}{g_q^-(a_1)} = \lim_{q \rightarrow 0}\mbf{E}_x[e^{-q\tau_1}] = 1$ if $x \ge a_1$, and  $\lim_{q \rightarrow 0}\frac{q_q^-(x)}{b_1^-} =\lim_{q \rightarrow 0}\mbf{E}_{a_1}[e^{-q\tau_x}] = 1$ if $x \le a_1$. Then by Lemma \ref{l0730b} and Theorem \ref{tpm} we have
 \beqnn
 \lim_{q \rightarrow 0}\mbf{P}_x\{X_{e_q} \in \d z\} = \lim_{q \rightarrow 0} \frac{q e^{-\delta_0^+(a_1 - z)}}{l_0\sigma_0^2 b_1^- (1-c^-_1)}  g_q^-(x) \d z  \frac{2e^{-2\mu_0(a_1 - z)/ \sigma_0^2}}{\sigma_0^2 \overline{F}_{1}}\d z.
 \eeqnn 
 Then the result follows.
\qed

	
\begin{prop}\label{p312b}
		If $\mu_0 > 0$ and $\mu_n < 0$, then for $x \in \mathbb{R}$, $z \in [a_i, a_{i+1}]$ with $i = 1, \cdots, n - 1$ and $F_i$ and $\overline{F}_{i + 1}$ given by \eqref{F} and \eqref{F1}, respectively,  we have
		\beqnn
		\lim_{q \rightarrow 0}\mbf{P}_x\{X_{e_q} \in \d z\} =  \frac{2e^{-2\mu_i(a_{i+1} - z)/\sigma_i^2}}{\sigma_i^2[F_{i}e^{-2\mu_i( a_{i+1}-a_{i})/\sigma_i^2}+ \overline{F}_{i + 1}]}\d z.
		\eeqnn
\end{prop}

\proof 	
By \eqref{Ci}, Lemmas \ref{l0730a} and \ref{l0730b}, one can check that
 \beqnn
\lim_{q \rightarrow 0}\frac{l_iC_i}{q} = \frac{F_{i}e^{-2\mu_i( a_{i+1}-a_{i})/\sigma_i^2}+ \overline{F}_{i + 1}}{2} + \frac{a_{i+1} - a_i}{ \sigma_i^2}{\bf 1}_{\{\mu_i = 0\}}.
\eeqnn 
For the case of $x \le z$, by Corollary \ref{exptau} we have 
\beqnn
\lim_{q \rightarrow 0}\frac{g_q^-(z)}{g_q^-(a_{i+1})} = \lim_{q\rightarrow 0}\mbf{E}_{a_{i+1}}[e^{-q\tau_z}] = 1
\eeqnn
and
\beqnn
\lim_{q \rightarrow 0}\frac{g_q^+(x)}{g_q^+(a_i)} = \begin{cases}
	\lim_{q \rightarrow 0}\mbf{E}_x[e^{-q\tau_i}] = 1, & x \le a_i,\\
	\lim_{q \rightarrow 0}1/\mbf{E}_{a_i}[e^{-q\tau_x}] = 1, & x > a_i.
\end{cases}
\eeqnn 
Then by Theorem \ref{tpm}, one obtains that
\beqnn
\lim_{q \rightarrow 0}\mbf{P}_x\{X_{e_q}\in \d z\} \ar=\ar \frac{qe^{-2\mu_i(a_{i+1} - z)/\sigma_i^2}	g_q^+(x)g_q^-(z)}{C_i l_i\sigma_i^2g_q^+(a_i)g_q^-(a_{i+1})} \d z \cr
\ar=\ar \frac{2e^{-2\mu_i(a_{i+1} - z)/\sigma_i^2}}{\sigma_i^2[F_{i}e^{-2\mu_i( a_{i+1}-a_{i})/\sigma_i^2}+ \overline{F}_{i + 1} + \frac{a_{i+1} - a_i}{\sigma_i^2}{\bf 1}_{\{\mu_i = 0\}}]}\d z.
\eeqnn 

 Similarly, for the case of $x > z$, one can check that $\lim_{q \rightarrow 0}\frac{g_q^-(x)g_q^+(z)}{g_q^+(a_i)g_q^-(a_{i+1})} = 1$. Then the result follows from the above and Theorem \ref{tpm}.
\qed

The following result can be obtained by the method in the proofs of Propositions \ref{p312a} and \ref{p312b}, which is stated without any proof.
\begin{prop}\label{p312}
	If $\mu_0 > 0$ and $\mu_n < 0$, then for $z \ge a_n$, $x \in \mathbb{R}$ and $F_n$  given by \eqref{F}, we have
	\beqnn
	\lim_{q \rightarrow 0}\mbf{P}_x\{X_{e_q} \in \d z\} = \frac{2e^{2\mu_n(z-a_n)/ {\sigma_n}^2}}{\sigma_n^2F_{n}}\d z.
	\eeqnn
\end{prop}

{\bf Proof of Corollary \ref{main1}}
  For any $i = 2, \cdots, n-1$, by the definitions of $F_i$ and $\overline{F}_{i + 1}$ in \eqref{F} and \eqref{F1}, we have
 \beqnn
	F_i \ar=\ar \overline{F}_i - \frac{2(a_{i+1} - a_i)}{\sigma_i^2}{\bf 1}_{\{\mu_i = 0\}} - e^{2\mu_{i }(a_{i+1} - a_{i })/\sigma_{i }^2}\overline{F}_{i + 1}\cr
	\ar\ar  + \frac{2(a_i - a_{i - 1})}{\sigma_{i - 1}^2}{\bf 1}_{\{\mu_{i - 1} = 0\}} + e^{-2\mu_{i - 1}(a_i - a_{i - 1})/\sigma_{i - 1}^2}F_{i - 1},
	\eeqnn 
	which implies that
	\beqnn
	L_i = e^{-2\mu_{i - 1}(a_i - a_{i - 1})/\sigma_{i - 1}^2}L_{i-1}
	\eeqnn 
	with 
	\beqlb\label{eqLi}
	L_i = F_i + e^{2\mu_{i }(a_{i+1} - a_{i })/\sigma_{i }^2}\overline{F}_{i + 1} +  \frac{2(a_{i+1} - a_i)}{\sigma_i^2}{\bf 1}_{\{\mu_i = 0\}}.
	\eeqlb 
	Then we have 
	\beqnn
	L_i = \exp\left\{-\sum_{\ell = 1}^{i-1}\frac{2\mu_\ell(a_{\ell + 1} - \ell)}{\sigma_\ell^2}\right\}L_1, \quad i = 1, \cdots, n-1
	\eeqnn 
	with $\sum_{\ell = 2}^1 = 0.$  Moreover,  
 	\beqnn
	L_1 = F_1 + e^{2\mu_{1 }(a_{2} - a_{1})/\sigma_{1 }^2}\overline{F}_{2} + \frac{2(a_2 - a_1)}{\sigma_i^1}{\bf 1}_{\{\mu_1 = 0\}} = \overline{F}_1
	\eeqnn  
	and
	\beqnn
	L_n:= F_n = \overline{F}_n + e^{-2\mu_{n - 1}(a_n - a_{n - 1})/\sigma_{n - 1}^2}F_{n - 1} = e^{-2\mu_{n - 1}(a_n - a_{n - 1})/\sigma_{n - 1}^2}L_{n-1}.
	\eeqnn 
    It follows that
	\beqlb\label{Li}
	L_i = \exp\left\{-\sum_{\ell = 1}^{i-1}\frac{2\mu_{\ell}(a_{\ell + 1} - a_{\ell})}{\sigma_{\ell}^2}\right\} \overline{F}_1, \quad i = 1, \cdots, n.
	\eeqlb 
 On the other hand, by the Terminal-Value Theorem (see, e.g., Schiff \cite[Theorem 2.36]{S99}) and \eqref{e_q}, one sees that
 \beqlb\label{eq0317}
 \lim_{t \rightarrow \infty}\mbf{P}_x\{X_t \in \d z\} = \lim_{q \rightarrow 0}\mbf{P}_x\{X_{e_q} \in \d z\}.
 \eeqlb
Then it follows from \eqref{eqLi}, \eqref{Li}, \eqref{eq0317} and Propositions \ref{p312a}-\ref{p312} that for any $x \in \mbb{R}$, we have
 \beqnn
 \ar\ar \lim_{t \rightarrow \infty}\mbf{P}_x\{X_t \in \d z\} \cr
 \ar\ar\quad= \begin{cases}
 	\frac{2e^{-2\mu_0(a_1 - z)/ \sigma_0^2}}{\sigma_0^2 \overline{F}_{1}}\d z, & z \le a_1,\\
 	\frac{2e^{2\mu_i(z - a_i)/\sigma_i^2}}{\sigma_i^2L_i}\d z,  & z \in (a_i, a_{i+1}], i = 1, \cdots, n - 1,\\
 	\frac{2e^{2\mu_n(z-a_n)/ {\sigma_n}^2}}{\sigma_n^2L_{n}}\d z, & z > a_n 
 \end{cases} \cr
\ar\ar\quad=  \frac{2}{\overline{F}_{1}}\begin{cases}
	\frac{1}{\sigma_0^2}e^{-2\mu_0(a_1 - z)/ \sigma_0^2}\d z, & z \le a_1,\\
	\frac{1}{\sigma_i^2}\exp\left\{\sum_{\ell = 1}^{i-1}\frac{2\mu_{\ell}(a_{\ell + 1} - a_{\ell})}{\sigma_{\ell}^2} + \frac{2\mu_i(z - a_i)}{\sigma_i^2}\right\}\d z, & z \in (a_i, a_{i+1}], i = 1, \cdots, n 
	\end{cases}
 \eeqnn
 with $a_{n+1} = +\infty.$ The proof is complete.
\qed

\section{Proof of Corollary \ref{t0227}}\label{SP}

In this section, we assume that $\mu_0 < 0$ and $\mu_n > 0$. Under this condition, $\mbf{P}_y\{\lim_{t \rightarrow \infty} X_t = -\infty\}$ is given in this section. The limit depends on the initial value $y$.

 \begin{lem}\label{p0226}
 Assuming that $\mu_n > 0$ and $\mu_i \neq 0$ for any $i = 1, \cdots, n$, we have
 	\beqlb\label{eq0307a}
 	\lim_{q \rightarrow 0}c_i^- = \begin{cases}
 		1 - \frac{\mu_i/\sigma_i^2}{\mu_{i-1}/\sigma_{i-1}^2}\frac{A_i}{A_i+B_i}, & \mu_{i - 1} > 0,\\
 		\frac{\mu_i/\sigma_i^2}{\mu_{i-1}/\sigma_{i-1}^2}\frac{A_i}{A_i+B_i}, & \mu_{i - 1} < 0,
 	\end{cases}
 	\eeqlb
  where $(A_i, B_i)_{i = 1, \cdots, n}$ are given by \eqref{ABi}.
 \end{lem} 

\proof By \eqref{delta}, one sees that 
\beqlb 
\lim_{q \rightarrow 0}\delta_i^- = \frac{|\mu_i| - \mu_i}{\sigma_i^2}, \quad \lim_{q \rightarrow 0}\delta_i^+ = \frac{|\mu_i| + \mu_i}{\sigma_i^2}, \quad \lim_{q \rightarrow 0}l_i = \frac{|\mu_i|}{\sigma_i^2}
\eeqlb
for any $i = 0, \cdots, n-1.$ Moreover,
\beqnn
\lim_{q \rightarrow 0} \frac{\delta_i^+ - \delta_{i+1}^+}{\delta_i^+ + \delta_i^-} =\frac12\left[1 + \frac{\mu_i}{|\mu_i|} - \frac{\sigma_i^2}{|\mu_i|}\frac{|\mu_{i+1}| + \mu_{i + 1}}{\sigma_{i+1}^2}\right], \quad \lim_{q \rightarrow 0}\frac{\delta_i^- + \delta_i^+}{\delta_{i-1}^+ + \delta_{i-1}^-} = \frac{|\mu_i|}{\sigma_i^2}\cdot\frac{\sigma_{i - 1}^2}{|\mu_{i - 1}|}.
\eeqnn
Recall that $\mu_n > 0.$ Then
\beqnn
\lim_{q \rightarrow 0}c_{n}^- \ar=\ar \lim_{q \rightarrow 0}\frac{\delta_{n-1}^+ - \delta_n^+}{\delta_{n-1}^+ + \delta_{n-1}^-} = \frac12\left[1 + \frac{\mu_{n-1}}{|\mu_{n-1}|} - \frac{\sigma_{n - 1}^2}{|\mu_{n-1}|}\frac{|\mu_{n}| + \mu_{n}}{\sigma_{n}^2}\right] \cr
\ar=\ar \frac12\left[1 + \frac{\mu_{n-1}}{|\mu_{n-1}|}\right] - \frac{\mu_n/\sigma_n^2}{|\mu_{n - 1}|/\sigma_{n - 1}^2}\cr
\ar=\ar \begin{cases}
	1 - \frac{\mu_n/\sigma_n^2}{\mu_{n - 1}/\sigma_{n - 1}^2}\frac{A_n}{A_n+B_n}, & \mu_{n-1} > 0,\\
	\frac{\mu_n/\sigma_n^2}{\mu_{n - 1}/\sigma_{n - 1}^2}\frac{A_n}{A_n + B_n}, & \mu_{n-1} < 0.
\end{cases}
\eeqnn
Then, \eqref{eq0307a} holds for $c_n^-$.  
Assuming that \eqref{eq0307a} holds for $c_i^-$, it suffices to show that it also holds for $i-1$. Notice that
\beqnn
\ar\ar \lim_{q \rightarrow 0}\frac{  c^-_{i}}{c^-_{i} +(1-c^-_i)e^{2l_{i-1}(a_{i} - a_{i-1})} } \cr
\ar\ar\quad =  \begin{cases}
	\frac{(1 - \frac{\mu_i/\sigma_i^2}{\mu_{i- 1}/\sigma_{i - 1}^2} \frac{A_i}{A_i+B_i})e^{-\mu_{i - 1}(a_i - a_{i - 1})/\sigma_{i - 1}^2} }{(1 - \frac{\mu_i/\sigma_i^2}{\mu_{i - 1}/\sigma_{i - 1}^2}\frac{A_i}{A_i+B_i})e^{-\mu_{i - 1}(a_i - a_{i - 1})/\sigma_{i - 1}^2} + \frac{\mu_i/\sigma_i^2}{\mu_{i - 1}/\sigma_{i - 1}^2}\frac{A_i}{A_i+B_i} e^{\mu_{i - 1}(a_i - a_{i - 1})/\sigma_{i - 1}^2}}, & \mu_{i-1} > 0,\\
	\frac{ \frac{\mu_i/\sigma_i^2}{\mu_{i - 1}/\sigma_{i - 1}^2}\frac{A_i}{A_i+B_i}e^{\mu_{i - 1}(a_i - a_{i - 1})/\sigma_{i - 1}^2} }{ \frac{\mu_i/\sigma_i^2}{\mu_{i - 1}/\sigma_{i - 1}^2}\frac{A_i}{A_i+B_i}e^{\mu_{i - 1}(a_i - a_{i - 1})/\sigma_{i - 1}^2} + (1 -\frac{\mu_i/\sigma_i^2}{\mu_{i - 1}/\sigma_{i - 1}^2}\frac{A_i}{A_i+B_i}) e^{-\mu_{i - 1}(a_i - a_{i - 1})/\sigma_{i - 1}^2}}, & \mu_{i-1} < 0
\end{cases}\cr
\ar\ar\quad =  \begin{cases}
	\frac{B_{i-1}}{A_{i-1} + B_{i-1}}, & \mu_{i-1} > 0,\\
	\frac{A_{i-1}}{A_{i-1} + B_{i-1}}, & \mu_{i-1} < 0.
\end{cases} 
\eeqnn 
It follows that
\beqnn    
\lim_{q \rightarrow 0} c_{i-1}^- \ar=\ar \lim_{q \rightarrow 0}\frac{1}{2l_{i-2}}\left[(\delta_{i-2}^+ - \delta_{i-1}^+) +  \frac{ 2l_{i-1} c^-_{i}}{c^-_{i} +(1-c^-_i)e^{2l_{i-1}(a_{i} - a_{i-1})} }\right]\cr
\ar=\ar \frac{\sigma_{i-2}^2}{2|\mu_{i-2}|}\left[\frac{|\mu_{i-2}| + \mu_{i-2}}{\sigma_{i-2}^2} - \frac{|\mu_{i-1}| + \mu_{i-1}}{\sigma_{i-1}^2} + \frac{2|\mu_{i-1}|}{\sigma_{i-1}^2}\lim_{q \rightarrow 0}\frac{  c^-_{i}}{c^-_{i} +(1-c^-_i)e^{2l_{i-1}(a_{i} - a_{i-1})} }\right]\cr
\ar=\ar  \begin{cases}
	1- \frac{\mu_{i-1}/\sigma_{i-1}^2}{\mu_{i-2}/\sigma_{i-2}^2}\frac{ A_{i-1}}{A_{i-1} + B_{i-1}}, & \mu_{i-2} > 0,\\
	\frac{\mu_{i-1}/\sigma_{i-1}^2}{\mu_{i-2}/\sigma_{i-2}^2} \frac{A_{i-1}}{A_{i-1} + B_{i-1}}, & \mu_{i-2} < 0.
\end{cases}
\eeqnn
Then the proof is completed by induction method. 
\qed 
 
\begin{prop}\label{p0806}
		Assume that $\mu_n > 0$. Then we have
		\beqnn
		\lim_{q \rightarrow 0}c_i^- = \begin{cases}
			1 - \frac{\mu_i/\sigma_i^2{\bf 1}_{\{\mu_i \neq 0\}} + {\bf 1}_{\{\mu_i = 0\}}}{\mu_{i-1}/\sigma_{i-1}^2}\frac{A_i}{A_i + B_i}, & \mu_{i-1} > 0,\\
			\frac{\mu_i/\sigma_i^2{\bf 1}_{\{\mu_i \neq 0\}} + {\bf 1}_{\{\mu_i = 0\}}}{\mu_{i-1}/\sigma_{i-1}^2}\frac{A_i}{A_i + B_i}, & \mu_{i - 1}< 0,
		\end{cases}
		\eeqnn 
		  where $(A_i, B_i)_{i = 1, \cdots, n}$ are given by \eqref{ABi}.
\end{prop}

\proof
Assume that $\mu_{n - 1} = 0.$ Then we have
\beqnn
\lim_{q \rightarrow 0}c_n^- = \lim_{q \rightarrow 0}\frac{\delta_{n-1}^+ - \delta_n^+}{\delta_{n-1}^+ + \delta_{n-1}^-} = \frac{1}{2} - \lim_{q \rightarrow 0}\frac{\delta_n^+}{2l_{n-1}} = -\infty
\eeqnn
and $\lim_{q \rightarrow 0}l_{n-1}c_n^- = -\frac{\mu_n}{\sigma_n^2}$. When $\mu_{n - 2} > 0$, it follows that
\beqnn
\lim_{q \rightarrow 0} c_{n-1}^- \ar=\ar \lim_{q \rightarrow 0}\left[\frac{\delta_{n-2}^+ - \delta_{n-1}^+}{\delta_{n-2}^+ + \delta_{n-2}^-} + \frac{\delta_{n-1}^- + \delta_{n-1}^+}{\delta_{n-2}^+ + \delta_{n-2}^-}\frac{ c^-_{n}}{c^-_{n} + (1- c^-_{n})e^{(\delta_{n-1}^+ + \delta_{n-1}^-)(a_{n} - a_{n-1})}}\right] \cr
\ar=\ar 1 + \lim_{q \rightarrow 0}\frac{1}{2l_{n-2}}\frac{l_{n-1}c_n^-(1 + e^{2l_{n-1}(a_{n} - a_{n-1})}) - l_{n-1}e^{2l_{n-1}(a_{n} - a_{n-1})}}{l_{n-1}c^-_{n}\frac{1 - e^{2l_{n-1}(a_{n} - a_{n-1})}}{l_{n-1}} + e^{2l_{n-1}(a_{n} - a_{n-1})}}\cr
\ar=\ar 1 - \frac{\mu_n/\sigma_n^2}{\mu_{n-2}/\sigma_{n-2}^2} \frac{1}{1 + 2\frac{\mu_n}{\sigma_n^2}(a_n - a_{n-1})}\cr
\ar=\ar 1 - \frac{1}{\mu_{n-2}/\sigma_{n-2}^2}\frac{A_{n-1}}{A_{n-1} + B_{n-1}}.
\eeqnn 
For $\mu_{n-2} < 0$, we have
\beqnn
\lim_{q \rightarrow 0} c_{n-1}^- \ar=\ar \lim_{q \rightarrow 0}\left[\frac{\delta_{n-2}^+ - \delta_{n-1}^+}{\delta_{n-2}^+ + \delta_{n-2}^-} + \frac{\delta_{n-1}^- + \delta_{n-1}^+}{\delta_{n-2}^+ + \delta_{n-2}^-}\frac{ c^-_{n}}{c^-_{n} + (1- c^-_{n})e^{(\delta_{n-1}^+ + \delta_{n-1}^-)(a_{n} - a_{n-1})}}\right] \cr
\ar=\ar \lim_{q \rightarrow 0}\frac{1}{2l_{n-2}}\frac{l_{n-1}c_n^-(1 + e^{2l_{n-1}(a_{n} - a_{n-1})}) - l_{n-1}e^{2l_{n-1}(a_{n} - a_{n-1})}}{l_{n-1}c^-_{n}\frac{1 - e^{2l_{n-1}(a_{n} - a_{n-1})}}{l_{n-1}} + e^{2l_{n-1}(a_{n} - a_{n-1})}}\cr
\ar=\ar \frac{\mu_n/\sigma_n^2}{\mu_{n-2}/\sigma_{n-2}^2} \frac{1}{1 + 2\frac{\mu_n}{\sigma_n^2}(a_n - a_{n-1})}\cr
\ar=\ar \frac{1}{\mu_{n-2}/\sigma_{n-2}^2}\frac{A_{n-1}}{A_{n-1} + B_{n-1}}.
\eeqnn 
Then the result holds for $i = n-1.$ 

 For the case of $\mu_i \neq 0$ and $\mu_{i-1} = 0$, we have
\beqlb\label{eq0806}
\lim_{q \rightarrow 0}l_{i-1}c_i^- \ar=\ar \frac{1}{2}\lim_{q \rightarrow 0}\frac{c_{i+1}^-(\delta_i^- + \delta_i^+ e^{2l_i(a_{i+1} - a_i)}) - \delta_i^+ e^{2l_i(a_{i+1} - a_i)}}{c^-_{i+1}(1 - e^{2l_i(a_{i+1} - a_i)}) +  e^{2l_i(a_{i+1} - a_i)}}\cr
\ar=\ar 
-\frac{\frac{\mu_{i+1}}{\sigma_{i+1}^2}A_{i+1}e^{2\mu_i(a_{i+1} - a_i)/\sigma_{i}^2}}{A_{i+1} + B_{i+1} + \frac{\mu_{i+1}/\sigma_{i+1}^2}{\mu_i/\sigma_i^2}A_{i+1}(e^{2\mu_i(a_{i+1} - a_i)/\sigma_i^2} - 1)}\cr
\ar=\ar - \frac{\mu_i}{\sigma_i^2}\frac{A_i}{A_i + B_i}.
\eeqlb
 For $\mu_{i - 2} > 0$, we have
\beqnn
\lim_{q \rightarrow 0}c_{i - 1}^- \ar=\ar  \lim_{q \rightarrow 0}  \left[\frac{\delta_{i-2}^+ - \delta_{i-1}^+}{\delta_{i-2}^+ + \delta_{i-2}^-} + \frac{\delta_{i-1}^- + \delta_{i-1}^+}{\delta_{i-2}^+ + \delta_{i-2}^-}\frac{ c^-_{i}}{c^-_{i} + (1- c^-_{i})e^{(\delta_{i-1}^+ + \delta_{i-1}^-)(a_{i} - a_{i-1})}}\right] \cr
\ar=\ar 1 + \lim_{q \rightarrow 0} \frac{l_{i - 1}}{2l_{i - 2}}\frac{c_i^-(1 + e^{2l_{i - 1}(a_i - a_{i - 1})})  - e^{2l_{i - 1}(a_i - a_{i - 1})}}{c^-_{i} (1 - e^{2l_{i-1}(a_{i} - a_{i-1})}) + e^{2l_{i-1}(a_{i} - a_{i-1})}}\cr
\ar=\ar  1 + \lim_{q \rightarrow 0} \frac{1}{2l_{i - 2}}\frac{l_{i - 1} c_i^-(1 + e^{2l_{i - 1}(a_i - a_{i - 1})})  - l_{i - 1} e^{2l_{i - 1}(a_i - a_{i - 1})}}{l_{i - 1}c^-_{i} \frac{1 - e^{2l_{i-1}(a_{i} - a_{i-1})}}{l_{i - 1}} + e^{2l_{i-1}(a_{i} - a_{i-1})}}\cr
\ar=\ar 1 + \frac{1}{\mu_{i - 2}/\sigma_{i-2}^2}\frac{\lim_{q \rightarrow 0} l_{i - 1}c_i^-}{1 - 2(a_i - a_{i - 1})\lim_{q \rightarrow 0}l_{i - 1}c_i^-}\cr
\ar=\ar 1 -   \frac{\mu_{i}/\sigma_{i}^2}{\mu_{i - 2}/\sigma_{i-2}^2}\frac{A_{i}}{A_{i} + B_{i} + 2\frac{\mu_{i}}{\sigma_{i}^2}(a_i - a_{i - 1})A_{i}}\cr
\ar=\ar 1 - \frac{1}{\mu_{i-2}/\sigma_{i-2}^2}\frac{A_{i-1}}{A_{i-1} + B_{i-1}}
\eeqnn 
For the case of $\mu_{i - 2} < 0,$ one obtains that
\beqnn
\lim_{q \rightarrow 0}c_{i - 1}^- \ar=\ar  \lim_{q \rightarrow 0}  \left[\frac{\delta_{i-2}^+ - \delta_{i-1}^+}{\delta_{i-2}^+ + \delta_{i-2}^-} + \frac{\delta_{i-1}^- + \delta_{i-1}^+}{\delta_{i-2}^+ + \delta_{i-2}^-}\frac{ c^-_{i}}{c^-_{i} + (1- c^-_{i})e^{(\delta_{i-1}^+ + \delta_{i-1}^-)(a_{i} - a_{i-1})}}\right] \cr
\ar=\ar \lim_{q \rightarrow 0} \frac{l_{i - 1}}{2l_{i - 2}}\frac{c_i^-(1 + e^{2l_{i - 1}(a_i - a_{i - 1})})  - e^{2l_{i - 1}(a_i - a_{i - 1})}}{c^-_{i} (1 - e^{2l_{i-1}(a_{i} - a_{i-1})}) + e^{2l_{i-1}(a_{i} - a_{i-1})}}\cr
\ar=\ar -\frac{1}{\mu_{i-2}/\sigma_{i-2}^2} \frac{\lim_{q \rightarrow 0}l_{i-1}c_i^-}{\lim_{q \rightarrow 0}l_{i-1}c_i^-\lim_{q \rightarrow 0}\frac{1 - e^{2l_{i-1}(a_i - a_{i - 1})}}{l_{i-1}} + 1}\cr
\ar=\ar  \frac{\mu_{i}/\sigma_{i}^2}{\mu_{i - 2}/\sigma_{i-2}^2}\frac{A_{i}}{A_{i} + B_{i} + 2\frac{\mu_{i}}{\sigma_{i}^2}(a_i - a_{i - 1})A_{i}}\cr
\ar=\ar \frac{1}{\mu_{i-2}/\sigma_{i-2}^2}\frac{A_{i-1}}{A_{i-1} + B_{i-1}}.
\eeqnn 
The result follows from induction method and Lemma \ref{p0226}.
\qed 

\begin{cor}\label{c0806}
	Assume that $\mu_n > 0$ and $\mu_0 < 0$. Then we have
	\beqnn
	\lim_{q \rightarrow 0}c_1^- = \frac{\mu_1/\sigma_1^2{\bf 1}_{\{\mu_1 \neq 0\}} + {\bf 1}_{\{\mu_1 = 0\}}}{\mu_{0}/\sigma_{0}^2}\frac{A_1}{A_1 + B_1},
	\eeqnn 
		where $(A_1, B_1)$ are given by \eqref{ABi}.
\end{cor}

\begin{prop}\label{t0226}
	Assume that $\mu_0 < 0$ and $\mu_n > 0.$ Then
	\beqnn
	\lim_{x \rightarrow -\infty}\mbf{P}_{a_1}\{\tau_x < \infty\} = \frac{A_1 + B_1}{ \left(1 - \frac{\mu_1/\sigma_1^2{\bf 1}_{\{\mu_1 \neq 0\}} + {\bf 1}_{\{\mu_1 = 0\}}}{\mu_0/\sigma_0^2}\right)A_1 + B_1},
	\eeqnn 
	where $(A_1, B_1)$ are given by \eqref{ABi}.
\end{prop} 
\proof By \eqref{delta}, we have 
\beqlb\label{eq0226}
\lim_{q \rightarrow 0} \delta_0^- = \frac{|\mu_0| - \mu_0}{\sigma_0^2} = -\frac{2\mu_0}{\sigma_0^2}, \qquad \lim_{q \rightarrow 0} \delta_0^+ = \frac{|\mu_0| + \mu_0}{\sigma_0^2} = 0
\eeqlb
and
\beqlb\label{eq0226a}
\lim_{q \rightarrow 0} \delta_n^- = \frac{|\mu_n| - \mu_n}{\sigma_n^2} = 0, \qquad \lim_{q \rightarrow 0} \delta_n^+ = \frac{|\mu_n| + \mu_n}{\sigma_n^2} = \frac{2\mu_n}{\sigma_n^2}.
\eeqlb
By Lemma \ref{l0809a}, for any $x  < a_1$, we have
\beqnn
\mbf{E}_{a_1} [e^{-q\tau_x}] = \frac{1}{(1 - c_1^-)e^{\delta_0^+(a_1 - x)} + c_1^- e^{-\delta_0^-(a_1 - x)}}.
\eeqnn
Then we have
\beqnn
\mbf{P}_{a_1}\{\tau_x < \infty\} \ar=\ar \lim_{q \rightarrow 0}\mbf{E}_{a_1}[e^{-q\tau_x}] = \lim_{q \rightarrow 0} \frac{1}{(1 - c_1^-)e^{\delta_0^+(a_1 - x)} + c_1^- e^{-\delta_0^-(a_1 - x)}}\cr
\ar=\ar  \frac{1}{\lim_{q \rightarrow 0}(1 - c_1^-) + \lim_{q \rightarrow 0}c_1^- e^{2\mu_0(a_1 - x)/\sigma_0^2}},
\eeqnn 
which implies that
\beqnn
\lim_{x \rightarrow -\infty} \mbf{P}_{a_1}\{\tau_x < \infty\} = \frac{1}{ \lim_{q \rightarrow 0}(1 - c_1^-) }. 
\eeqnn 
The result follows from Corollary \ref{c0806}. \qed
 
{\bf Proof of Corollary \ref{t0227}}
 (i) For the case of $a_1 \ge y > x$, by the strong Markov property, we have
\beqnn
\mbf{E}_y[e^{-q\tau_x}] \ar=\ar \mbf{E}_y[e^{-q\tau_x}{\bf 1}_{\{\tau_x < \tau_1\}}] + \mbf{E}_y[e^{-q(\tau_x - \tau_1 + \tau_1)}{\bf 1}_{\{\tau_x \ge \tau_1\}}]\cr
\ar=\ar \mbf{E}_y[e^{-q\tau_x}{\bf 1}_{\{\tau_x < \tau_1\}}] + \mbf{E}_y[e^{-q\tau_1}{\bf 1}_{\{\tau_x \ge \tau_1\}}]\mbf{E}_{a_1}[e^{-q\tau_x}].
\eeqnn 
By \eqref{star4} and \eqref{star5}, one obtains that
\beqlb\label{eq0504a}
\lim_{q \rightarrow 0} \mbf{E}_y[e^{-q\tau_x}{\bf 1}_{\{\tau_x < \tau_1\}}] \ar=\ar  \lim_{q \rightarrow 0} \mbf{E}_y[e^{-q(\tau_x \wedge \tau_1)}{\bf 1}_{\{\tau_x < \tau_1\}}]\cr
\ar=\ar \lim_{q \rightarrow 0}e^{\mu_0(x - y)/\sigma^2}\frac{\sinh((a_1 - y)l_0)}{\sinh((a_1 - x)l_0)}\cr
\ar=\ar e^{\mu_0(x - y)/\sigma^2}\frac{\sinh(-(a_1 - y)\mu_0/\sigma_0^2)}{\sinh(-(a_1 - x)\mu_0/\sigma_0^2)}
\eeqlb
and
\beqlb\label{eq0504b}
\lim_{q \rightarrow 0}\mbf{E}_y[e^{-q\tau_1}{\bf 1}_{\{\tau_x \ge \tau_1\}}] \ar=\ar \lim_{q \rightarrow 0}\mbf{E}_y[e^{-q(\tau_1\wedge \tau_x)}{\bf 1}_{\{\tau_x \ge \tau_1\}}] \cr
\ar=\ar \lim_{q \rightarrow 0}e^{\mu_0(a_1 - y)/\sigma_0^2}\frac{\sinh((y-x)l_0)}{\sinh((a_1 - x)l_0)}\cr
\ar=\ar e^{\mu_0(a_1 - y)/\sigma_0^2}\frac{\sinh(-\mu_0(y-x)/\sigma_0^2)}{\sinh(-\mu_0(a_1 - x)/\sigma_0^2)}.
\eeqlb
Then
\beqnn
\mbf{P}_y\{\tau_x < \infty\} \ar=\ar \lim_{q \rightarrow 0}\mbf{E}_y[e^{-q\tau_x}]\cr
\ar=\ar e^{\mu_0(x - y)/\sigma^2}\frac{\sinh(-(a_1 - y)\mu_0/\sigma_0^2)}{\sinh(-(a_1 - x)\mu_0/\sigma_0^2)}\cr
\ar\ar + e^{\mu_0(a_1 - y)/\sigma_0^2}\frac{\sinh(-\mu_0(y-x)/\sigma_0^2)}{\sinh(-\mu_0(a_1 - x)/\sigma_0^2)} \mbf{P}_{a_1}\{\tau_x < \infty\}.
\eeqnn
Letting $x \rightarrow -\infty$ in the above, it follows from Theorem \ref{t0226} that
\beqnn
\lim_{x \rightarrow -\infty} \mbf{P}_y\{\tau_x < \infty\} \ar=\ar 1 - e^{2\mu_0(a_1 - y)/\sigma_0^2}+ e^{2\mu_0(a_1 - y)/\sigma_0^2}\lim_{x \rightarrow -\infty} \mbf{P}_{a_1}\{\tau_x < \infty\}\cr
\ar=\ar  1 - e^{2\mu_0(a_1 - y)/\sigma_0^2}+ \frac{e^{2\mu_0(a_1 - y)/\sigma_0^2}(A_1 + B_1)}{\left(1 - \frac{\mu_1/\sigma_1^2{\bf 1}_{\{\mu_1 \neq 0\}} + {\bf 1}_{\{\mu_1 = 0\}}}{\mu_0/\sigma_0^2}\right)A_1 + B_1}\cr
\ar=\ar \frac{A_1(y) + B_1(y)}{\left(1 - \frac{\mu_1/\sigma_1^2 {\bf 1}_{\{\mu_1 \neq 0\}} + {\bf 1}_{\{\mu_1 = 0\}}}{\mu_0/\sigma_0^2}\right)A_1 + B_1}     e^{\mu_0(a_1 - y)/\sigma_0^2}.
\eeqnn 
 
(ii) For the case of $y \in (a_{i-1}, a_{i}]$ with $i = 2, \cdots, n$ and $x < a_1$, by the strong Markov property, we have
\beqnn
\mbf{E}_y[e^{-q\tau_x}] = \mbf{E}_y[e^{-q(\tau_x - \tau_1 + \tau_1)}] = \mbf{E}_y[e^{-q\tau_1}]\mbf{E}_{a_1}[e^{-q \tau_x}].
\eeqnn
Letting $q \rightarrow 0$, we have
\beqlb\label{eq0304b}
\mbf{P}_y\{\tau_x < \infty\} = \lim_{q \rightarrow 0}\mbf{E}_y[e^{-q\tau_x}] = \lim_{q \rightarrow 0}\mbf{E}_y[e^{-q\tau_1}] \mbf{P}_{a_1}\{\tau_x < \infty\}. 
\eeqlb
It suffices to calculate $\lim_{q \rightarrow 0}\mbf{E}_y[e^{-q\tau_i}].$ 
By \eqref{eq0806}, for $y \in (a_{i-1}, a_i]$ we have
\beqlb\label{eq0304}
\ar\ar \lim_{q \rightarrow 0}c^-_{i}e^{\delta_{i-1}^-(y-a_i) }+(1-c^-_{i})e^{-\delta_{i-1}^+(y-a_i)}\cr
\ar\ar\quad = \lim_{q \rightarrow 0}l_{i-1}c_i^-\frac{e^{2l_{i-1}(y-a_i) } - 1}{l_{i-1}}e^{-\delta_{i-1}^+(y-a_i)} + e^{-\delta_{i-1}^+(y-a_i)}\cr
\ar\ar\quad =  1 + \frac{2\mu_i/\sigma_i^2(a_i - y)A_i}{A_i + B_i}\cr
\ar\ar\quad  =  \frac{A_{i}(y) + B_{i}(y)}{A_i + B_i}e^{\mu_{i-1}(a_i - y)/\sigma_{i-1}^2} 
\eeqlb
when $\mu_{i-1} = 0$. And by Proposition \ref{p0806}, one sees that
\beqnn
\ar\ar \lim_{q \rightarrow 0}c^-_{i}e^{\delta_{i-1}^-(y-a_i) }+(1-c^-_{i})e^{-\delta_{i-1}^+(y-a_i)}\cr
\ar\ar\quad = 1 + \frac{\mu_i/\sigma_i^2 {\bf 1}_{\{\mu_i \neq 0\}} + {\bf 1}_{\{\mu_i = 0\}}}{\mu_{i-1}/\sigma_{i-1}^2}\frac{A_i}{A_i+B_i}( e^{2\mu_{i-1}(a_i - y)/\sigma_{i-1}^2} - 1)\cr
\ar\ar\quad = \frac{A_{i}(y) + B_{i}(y)}{A_i + B_i}e^{\mu_{i-1}(a_i - y)/\sigma_{i-1}^2} 
\eeqnn 
when $\mu_{i - 1} \neq 0.$ Similar to above, by \eqref{eq0806} and Proposition \ref{p0806}, one can check that
\beqnn
 \lim_{q \rightarrow 0} \left[c^-_{\ell+1}e^{-\delta_\ell^-(a_{\ell+1} - a_\ell)} + (1- c^-_{\ell+1})e^{\delta_\ell^+(a_{\ell+1} - a_\ell)}\right] = \frac{A_\ell + B_\ell}{A_{\ell + 1} + B_{\ell + 1}}e^{-\mu_\ell(a_{\ell + 1} - a_\ell)/\sigma_\ell^2}.
\eeqnn 
Then it follows from the above and \eqref{b^-} that
\beqlb\label{eq0304a}
\lim_{q \rightarrow 0}\frac{b_i^-}{b_1^-} \ar=\ar \lim_{q \rightarrow 0}\prod_{\ell = 1}^{i-1}\frac{b_{\ell+1}^-}{b_\ell^-} = \lim_{q \rightarrow 0}\prod_{\ell = 1}^{i-1} \frac{1}{c^-_{\ell+1}e^{-\delta_\ell^-(a_{\ell+1} - a_\ell)} + (1- c^-_{\ell+1})e^{\delta_\ell^+(a_{\ell+1} - a_\ell)}}\cr
\ar=\ar \frac{A_{i} + B_{i}}{A_1 + B_1} e^{-\sum_{\ell = 1}^{i - 1}\mu_\ell (a_{\ell + 1} - a_\ell)/\sigma_\ell^2}.
\eeqlb
By \eqref{eq0304} and \eqref{eq0304a}, we have
\beqnn
\lim_{q \rightarrow 0}\mbf{E}_y[e^{-q\tau_1}] \ar=\ar \lim_{q \rightarrow 0}\frac{g_q^-(y)}{g_q^-(a_1)} = \lim_{q \rightarrow 0}\frac{b_i^-}{b_1^-}\bigg( c^-_{i}e^{\delta_{i-1}^-(y-a_i) }+(1-c^-_{i})e^{-\delta_{i-1}^+(y-a_i)}\bigg)\cr
\ar=\ar \frac{ A_{i}(y) + B_{i}(y)}{A_1 + B_1}e^{\mu_{i-1}(a_i - y)/\sigma_{i-1}^2-\sum_{\ell = 1}^{i-1}\mu_\ell(a_{\ell + 1} - a_\ell)/\sigma_\ell^2 }.
\eeqnn
It follows from the above, \eqref{eq0304b} and Proposition \ref{t0226} that
\beqnn
\lim_{x \rightarrow -\infty} \mbf{P}_y\{\tau_x < \infty\} \ar=\ar \lim_{q \rightarrow 0} \mbf{E}_y[e^{-q\tau_1}] \lim_{x \rightarrow -\infty} \mbf{P}_{a_1}\{\tau_x < \infty\}\cr
\ar=\ar \frac{ A_{i}(y) + B_{i}(y)}{\left(1 - \frac{\mu_1/\sigma_1^2}{\mu_0/\sigma_0^2}\right)A_1 + B_1}e^{\mu_{i-1}(a_i - y)/\sigma_{i-1}^2-\sum_{\ell = 1}^{i-1}\mu_\ell(a_{\ell + 1} - a_\ell)/\sigma_\ell^2 }.
\eeqnn 

(iii) For the case of $y > a_n$, \eqref{eq0304b} still holds.
Moreover, by \eqref{g_q^-} and Corollary \ref{exptau},
\beqnn
\lim_{q \rightarrow 0}\mbf{E}_y[e^{-q\tau_1}] \ar=\ar \lim_{q \rightarrow 0}\frac{g_q^-(y)}{g_q^-(a_1)} = \lim_{q \rightarrow 0}\frac{e^{-\delta_n^+(y - a_n)}}{b_1^-} = \lim_{q \rightarrow 0} e^{-\delta_n^+(y - a_n)}\prod_{\ell = 1}^{n-1}\frac{b_{\ell+1}^-}{b_\ell^-}\cr
\ar=\ar \frac{e^{-2\mu_n(y - a_n)/\sigma_n^2}}{\prod_{\ell = 1}^{n - 1}\frac{A_\ell + B_\ell}{A_{\ell + 1} + B_{\ell + 1}}e^{\mu_\ell(a_{\ell + 1} - a_\ell)/\sigma_\ell^2}}\cr
\ar=\ar \frac{1}{A_1 + B_1}e^{-\sum_{\ell = 1}^{n-1}\mu_\ell(a_{\ell + 1} - a_\ell)/\sigma_\ell^2-2\mu_n(y - a_n)/\sigma_n^2}. 
\eeqnn
Then by the above and Proposition \ref{t0226}, we have
\beqnn
\lim_{x \rightarrow -\infty}\mbf{P}_y\{\tau_x < \infty\} \ar=\ar \lim_{q \rightarrow 0}\mbf{E}_y[e^{-q\tau_1}] \lim_{x \rightarrow -\infty} \mbf{P}_{a_1}\{\tau_x < \infty\}\cr
\ar=\ar \frac{1}{\left(1 - \frac{\mu_1/\sigma_1^2}{\mu_0/\sigma_0^2}\right)A_1 + B_1}e^{-\sum_{\ell = 1}^{n-1}\mu_\ell(a_{\ell + 1} - a_\ell)/\sigma_\ell^2-2\mu_n(y - a_n)/\sigma_n^2}. 
\eeqnn 
Notice that 
\beqnn
\mbf{P}_y\{\lim_{t \rightarrow \infty}X_t = -\infty\} = \lim_{x \rightarrow -\infty}\mbf{P}_y\{\tau_x < \infty\}.
\eeqnn
Then the result follows.
\qed

\end{document}